\documentclass[journal,onecolumn,12pt]{IEEEtran} 

\usepackage[utf8]{inputenc}
\usepackage[T1]{fontenc}
\usepackage[hidelinks]{hyperref}
\usepackage{url}
\usepackage{ifthen}
\usepackage{cite}
\usepackage[cmex10]{amsmath} 
\usepackage{amsfonts,amssymb}
\usepackage{color}
\usepackage{ifthen}
\usepackage[arxiv]{optional} 

\usepackage{cite}
\usepackage{amsmath,amssymb,amsfonts}
\usepackage{algorithmic}
\usepackage{graphicx}
\usepackage{textcomp}
\usepackage{xcolor}
\usepackage{caption}
\usepackage{subcaption}
\usepackage{bbm}
\usepackage{graphics} 
\usepackage{epsfig} 
\usepackage{times} 
\usepackage{amsmath} 
\usepackage{amssymb}  
\usepackage{caption}
\usepackage{hyperref}
\usepackage{bbold}

\newcommand{\eq}[1]{\begin{align}#1\end{align}}
\newcommand{\seq}[1]{\begin{subequations}#1\end{subequations}}

\newcommand{\lb}[1]{\left\{ \begin{array}{ll} #1 \end{array} \right.}

\newcommand{\bit}[1]{\begin{itemize}#1\end{itemize}}
\newcommand{\E}{\mathbb{E}}
\newcommand{\p}{\mathbb{P}}

\newcommand{\cT}{[T]}

\newcommand{\cP}{\mathcal{P}}

\newcommand{\cX}{\mathcal{X}}

\newcommand{\cN}{\mathcal{N}}

\newcommand{\cA}{\mathcal{A}}

\newcommand{\tsigma}{\tilde{\sigma}}
\newcommand{\tgamma}{\tilde{\gamma}}

\newcommand{\defeq}{\buildrel\triangle\over =}

\newcommand{\nn}{\nonumber}

\DeclareMathAlphabet{\mathcal}{OMS}{cmsy}{m}{n}

\usepackage{lipsum}
\usepackage{amsfonts}
\usepackage{graphicx}
\usepackage{epstopdf}
\usepackage{algorithmic}
\ifpdf
  \DeclareGraphicsExtensions{.eps,.pdf,.png,.jpg}
\else
  \DeclareGraphicsExtensions{.eps}
\fi

\newtheorem{lemma}{Lemma}

\newtheorem{theorem}{Theorem}

\newtheorem{proof}{Proof}
\newtheorem{definition}{Definition}
\newtheorem{claim}{Claim}

\title{Sequential decomposition of stochastic Stackelberg games}

\author{Deepanshu Vasal\thanks{Department of Electrical Engineering, Northwestern University,
  ({dvasal@umich.edu}, \url{https://sites.google.com/view/dvasal/home}).\newline
  Part of the paper has been accepted at American Control Conference, 2022.}}

\usepackage{amsopn}

\ifpdf
\hypersetup{
  pdftitle={Stochastic Stackelberg games},
  pdfauthor={D. Vasal}
}
\fi

\begin{document}

\maketitle

\begin{abstract}
In this paper, we consider a discrete-time stochastic Stackelberg game with a single leader and multiple followers. Both the followers and the leader together have conditionally independent private types, conditioned on action and previous state, that evolve as controlled Markov processes. The objective is to compute the stochastic Stackelberg equilibrium of the game where the leader commits to a dynamic strategy. Each follower's strategy is the best response to the leader's strategies and other followers' strategies while the each leader's strategy is optimum given the followers play the best response. In general, computing such equilibrium involves solving a fixed-point equation for the whole game. In this paper, we present a backward recursive algorithm that computes such strategies by solving smaller fixed-point equations for each time $t$. Based on this algorithm, we compute stochastic Stackelberg equilibrium of a security example and a dynamics information design example used in~\cite{El17} (beeps).
\end{abstract}
%
%

\section{Introduction}
In the past decade, Stackelberg games have been used extensively in the security of real world systems such as to protect ports, airports and wildlife~\cite{BaGaAm09,Pietal08,Kietal09,CoSa06}. A Bayesian Stackelberg game is played between two players: a leader and an follower. The follower has a private type that only she observes, however, the leader knows the prior on that state. The leader commits to a strategy that is observable to the follower. The follower then plays a best response to follower's strategy to maximize its utility. Knowing that the follower will play a best response, the leader commits to and plays a strategy that maximizes its utility. Such pair of strategies is called a Stackelberg equilibrium. It is known that such strategies can provide higher utility to the leader than obtained in a Nash equilibrium of the game. Such games have been used in the real world by security agencies such as the US Coast Guard, the federal Air Marshals Service, and the Los Angeles Airport Police~\cite{Ta11}. Similar algorithms are used in wildlife protection in Uganda and Malaysia~\cite{Yaetal14}. 

Most of the above real world applications of Stackelberg equilibrium are based on single-shot Bayesian game models. However, in many practical scenarios, the follower and leader interact periodically, and also have private information, thus reducing the applicability of such models.
Such games comes under the class of dynamic games of asymmetric information, where both the leader and the follower privately observe conditionally independent controlled Markov processes, but observe each others' actions publicly. The reason such games are hard is because in such games the beliefs that come up across the game at any time $t$ are dependent on the strategies of the players before this time. Thus there is no notion of ``state" that can decompose the problem across time, and effectively there is no notion of dynamic programming. The space of strategies grows double exponential in time making solving for equilibria for such problems impossible for any reasonable time duration.  Recently, there has been results on sequential decomposition of certain classes of games of asymmetric information~\cite{VaSiAn16arxiv,VaAn16allerton,VaAn16cdc}. 

In repeated Stackelberg security games, there have been other approaches to mitigate this issue. Mareki et.al. in~\cite{Maetal12} study a Bayesian repeated Stackelberg game where they assume leaders are myopic, thus significantly simplifying the analysis of finding the equilibrium. Balcan et al in~\cite{Baetal15} consider a learning theoretic approach to study a repeated Stackelberg game between follower and leader where they use regret analysis to learn follower's types, and show sub-linear regret for both complete and partial information models. Authors in~\cite{Ely-goalpost,Ely-sequential-info-design,Basu2017Dynamic,Hung-dynamic,Farhadi2018,Hamid-Allerton,Asu2020} study a dynamic Stackelberg game where there is a sender who observes a static state privately and has a commitment power, and with both long-term-optimizing principal and long-term-optimizing followers. Farhadi and Teneketzis in~\cite{FaTe20} consider a model where on top of all the assumptions in previous papers, the state is also dynamically evolving. To the best of our knowledge, \cite{FaTe20} is the only work that considers a special case of a truly dynamic Stackelberg game and in general finding Stackelberg equilibria of such games is an open question.  

In this paper, we show there indeed is a way to decompose general stochastic dynamic Stackelberg games across time, where both the leader and the follower have private Markovian states that evolve as conditionally independent Markov processes. We provide a dynamic programming like sequential decomposition algorithm to compute equilibria with fully rational, forward looking follower and leader. Our algorithm consists of a backward recursive step which, for each time t and a belief state on the current state, $\pi_t$, involves solving a fixed-point equation for the follower and an optimization problem for the leader. This reduces the complexity of finding Markovian equilibria of such games from double exponential to linear in time. Based on this algorithm, we study a security game where we numerically find its Stackelberg equilibria. To the best of our knowledge, this is the first paper that provides a general treatment to compute Stackelberg equilibria of Stochastic Stackelberg games with asymmetric information. Part of the paper (without any proofs) was published in~\cite{Va22}.

The paper is structured as follows. We present our model in Section~\ref{sec:Model}. We discuss background material and solution concept in Section~\ref{sec:Prelim}. In Section~\ref{sec:Result}, we present our main result of providing an algorithm to compute Markovian equilibrium strategies. In Section~\ref{sec:Infinite}, we discuss an infinite horizon version of the problem. 
We conclude in Section~\ref{sec:Conclusion}. All proofs are presented in the Appendices.

 \subsection{Notation}
We use uppercase letters for random variables and lowercase for their realizations. For any variable, subscripts represent time indices and superscripts represent player indices. We use notation $A_{t:t'}$ to represent the vector $(A_t, A_{t+1}, \ldots A_{t'})$ when $t'\geq t$ or an empty vector if $t'< t$. We remove superscripts or subscripts if we want to represent the whole vector, for example $ A_t$  represents $(A_t^1, \ldots, A_t^N) $. In a similar vein, for any collection of sets $(\cX^i)_{i \in \cN}$, we denote $\times_{i\in\cN} \cX^i$ by $\cX$. 
For any finite set $\mathcal{S}$, $\Delta(\mathcal{S})$ represents the space of probability measures on $\mathcal{S}$ and $|\mathcal{S}|$ represents its cardinality. We denote by $P^g$ (or $\E^g$) the probability measure generated by (or expectation with respect to) strategy profile $g$. We denote the set of real numbers by $\mathbb{R}$. For a probabilistic strategy profile of players $(\sigma_t^i)_{i\in \cN}$ where the probability of action $a_t^i$ conditioned on $a_{1:t-1},x_{1:t}^i$ is given by $\sigma_t^i(a_t^i|a_{1:t-1},x_{1:t}^i)$, we use the notation $\sigma_t^{-i}(a_t^{-i}|a_{1:t-1},x_{1:t}^{-i})$ to represent $\prod_{j\neq i} \sigma_t^j(a_t^j|a_{1:t-1},x_{1:t}^j)$.
All equalities/inequalities involving random variables are to be interpreted in the \emph{a.s.} sense.
For mappings with range function sets $f: \cA \rightarrow (\mathcal{B} \rightarrow \mathcal{C})$ we use square brackets $f[a] \in \mathcal{B} \rightarrow \mathcal{C}$ to denote the image of $a\in\mathcal{A}$ through $f$ and parentheses $f[a](b) \in \mathcal{C}$ to denote the image of $b\in\mathcal{B}$ through $f[a]$.
A controlled Markov process with state $X_t$, action $A_t$, and horizon $[T]$ is denoted by $(X_t,A_t)_{t\in\cT}$.

\section{Model}
\label{sec:Model}
We consider an incomplete information stochastic Stackelberg game over a time horizon $[T]\defeq$ $\{1, 2, \ldots T\}$ with simultaneous moves and perfect recall as follows. Suppose there is one Stackelberg leader and $M$ (Nash) followers. The leader and the followers have private types, $x_t^{l} \in \cX^{l}$ for the leader, $x_t^{m,j} \in \cX^m$, for $j$th follower where $x_t^{f,i},x_t^{m,j}$ evolve as a conditionally independent controlled Markov processes in the following way, where for the Stackelberg leader, $M$ (Nash) major followers. Let $ x_t^m = x_t^{m[1:M]}$, $a_t^m = a_t^{m[1:M]}$

\eq{
P(x_t^{l},x_t^{m[1:M]}|a^{l,m}_{1:t-1},x^{l,m}_{1:t-1}) 
&=   Q^{l}(x_t^{l}|a_{t-1}^{l,m},x_{t-1}^{l,m}) \prod_{j=1}^M Q^{m,j}(x_t^{m,j}|a_{t-1}^{l,m},x_{t-1}^{l,m}),
}
where $Q^{l},Q^{m,j}$ are known kernels defined below. The Stackelberg leader takes action $a_t^{l}\in \cA^{l}$ at time $t$ on observing $x_{1:t}^{l}$, a major follower $j$ takes action based on $x_{1:t}^{m,j}$  at time $t$ on observing $a^{l,m}_{1:t-1}$ and $x_{1:t}^{m,j}$.
Here, $a^{l,m}_{1:t-1}$ is common information among players, and $x_{1:t}^{l}, x_{1:t}^{m,j}$ are private information of the Stackelberg leader, major follower $j$, respectively. 
 At the end of interval $t$, Stackelberg leader receives an instantaneous reward $R_t^{l}(x_t^{l,m},a_t^{l,m})$, major follower $j$ receives an instantaneous reward $R_t^{m,j}(x_t^{l,m},a_t^{l,m})$. 

The sets $\cA^{l},\cA^{m,j}, \cX^{l},\cX^{m,j} $ are assumed to be finite. Let $\sigma^i = ( \sigma^i_t)_{t \in [T]}$ be a probabilistic strategy of a player $i\in\{l,m[1:M] \}$ where $\sigma^{l}_t : ((\cA^l\times\cA^m)^{t-1}\times(\cX^l)^t \to \mathcal{P}(\cA^{l})$, such that the leader plays action $A_t^{l}$ according to $ A_t^{l} \sim \sigma^{l}_t(\cdot|a_{1:t-1}^{l,m},x_{1:t}^{l})$, and the follower plays action $A_t^{m,j}$ according to $ A_t^{m,j} \sim \sigma^{m,j}_t(\cdot|a_{1:t-1}^{l,m},x_{1:t}^{m,j})$, Let $ \sigma \defeq(\sigma^i)_{i\in \{l,f\}}$ be a strategy profile of all players. Suppose players discount their rewards by a discount factor $\delta\leq 1$.

\section{Preliminaries}
\label{sec:Prelim}
In this section, we first present the definition of a Stackelberg Perfect Equilibrium (SPE) which we will use in this paper. We then discuss the common agent approach that we will utilize in deriving an algorithm for finding an SPE.

\subsection{Stackelberg perfect equilibrium }
\label{sec:PBSE}
In this paper, we will consider followers' Markovian equilibrium policy that only depends on her current states $x_t^{m,j}$ and common marginal beliefs $\pi_t^{l},\pi_t^{m,j}, j=1\ldots M$, where $\pi_t^{l}(x_t^{l}) = P^{\sigma^{l},\sigma^{m}}(x_t^{l}|a_{1:t-1}^{l,m}), \pi_t^{m,j}(x_t^{m,j}) = P^{\sigma^{l},\sigma^{m}}(x_t^{m,j}|a_{1:t-1}^{l,m})$ i.e. $\pi_t^{l}$ is the common belief on the Stackelberg leader's state given the common information $(a_{1:t-1}^{l,m})$. Let ${\underline{\pi_t}} = \{ \pi_t^{l}, \pi_t^{m,j}, \}_{ j=1\ldots K }$ Thus, at equilibrium the Stackelberg leader's strategy $a_t^{l}\sim \tsigma^{l}_t(\cdot|\underline{\pi_t},x_t^{l})$ and major follower $j$'s strategy $A_t^{m,j}\sim \tsigma^{m,j}_t(\cdot|\underline{\pi_t},x_t^{m,j})$.\footnote{Note, however, that for the purpose of equilibrium, we allow for deviations in the space of all possible strategies that may depend on the entire observation history.} 
 This specifies a minor follower's best response at time $t$ given the history of the mean-field state and its private type up to time $t$ and the leader's strategy from time $t$ on-wards.   
 Note that this mapping specifies a complete policy for the follower for all time $t$. 
Similarly we define the best response of major follower $j$ as
\eq{
&BR_t^{m,j}(\underline{\pi_t},a_{1:t-1}^{l,m},x_{1:t}^{m,j},\sigma_{t:T}^{l},\sigma_{t:T}^{m,-j}) \nn\\
&:=\arg\max_{ \sigma^{m,j}} \E^{\sigma_{t:T}^{l}\sigma_{t:T}^{m,j}\sigma_{t:T}^{m,-j},\underline{\pi_t}} \big\{ \sum_{n=t}^T \delta^{n-t}R_n^{m,j}(X_n^{l,m},A_n^{l,m}) |\underline{\pi_t},a_{1:t-1}^{l,m},x_{1:t}^{m,j}\big\}
}
\eq{ 
&BR^{m,j}(\sigma^l,\sigma^{m,-j}) :=\bigcap_t \bigcap_{a_{1:t-1}^{m,j}}\bigcap_{x_{1:t}^{m,j}}  BR_t^{m,j}(\underline{\pi_t},a_{1:t-1}^{l,m},x_{1:t}^{m,j},\sigma_{t:T}^{l},\sigma_{t:T}^{m,-j},\sigma_{t:T}^{m,j}).
}
With some abuse of notation, we will also say $\sigma_t^{m,j}\in BR_t^{m,j}(\sigma_{1:T}^{l,m})$ if there exists $\hat{\sigma}^{m,j} \in BR^{m,j}(\sigma_{1:T}^{l,m})$ such that $\sigma_t^{m,j} = \hat{\sigma}_t^{m,j}$. Similarly we say $\sigma_t^{m,j}\in BR_t^{m,j}(\sigma_{1:T}^{l},\sigma_{1:T}^{m,-j},\sigma_{1:T}^{m,j})$ if there exists $\hat{\sigma}^{m,j} \in BR^{m,j}(\sigma_{1:T}^{l},\sigma_{1:T}^{m,-j},\sigma_{1:T}^{m,j})$ such that $\sigma_t^{m,j} = \hat{\sigma}_t^{m,j}$

We now define best response of Stackelberg leader as follows
\eq{
BR^{l}(\sigma^{m}) &:=\bigcap_t \bigcap_{a_{1:t-1}^{l}}\bigcap_{x_{1:t}^{l}}  \arg\max_{ \sigma^{l}} \E^{\sigma^{l},\sigma^{m},\underline{\pi_t}} \big\{ \sum_{n=t}^T \delta^{n-t}R_n^{l}(X_n^{l,m},A_n^{l,m}) |\underline{\pi_t},a_{1:t-1}^{l,m},x_{1:t}^{l}\big\},\\
&\text{ where, } \hat{\sigma}^{m,j} \in BR^{m,j}(\sigma^{l},\hat{\sigma}^m),\hat{\sigma}^m\in BR^{m,j}(\sigma^{l},\hat{\sigma}^m)
}

\begin{definition}
\label{def:ME}
 $(\tsigma^l,\tsigma^{m})$ is a Stackelberg Perfect Equilibrium (SPE) if

[(a)]: $\tsigma^{l} \in BR^{l}(\tsigma^m)$

[(b)]:
$\tsigma^{m} \in BR^{n}(\tsigma^l,\tsigma^m)$

\end{definition}

\subsection{Common agent approach}
We recall that in general, the leader and the followers generate their actions at time $t$ as follows, $A_t^{l}\sim \sigma_t^{l}(\cdot|a_{1:t-1}^{l,m},x_{1:t}^{l})$, \\$A_t^{m,j}\sim \sigma_t^{m,j}(\cdot|a_{1:t-1}^{l,m},x_{1:t}^{m,j})$.
An alternative way to view the problem is as follows. As is done in the common information approach~\cite{NaMaTe13}, at time $t$, a fictitious common agent observes the common information $a_{1:t-1}^{l,m}$ and generates prescription functions $\gamma_t = (\gamma_t^{l},\gamma_t^m) = \psi_t[a_{1:t-1}^{l,m}]$, where $ \gamma^{m,j} = \{\gamma_t^m\}_{j=1\ldots M}$. The Stackelberg leader uses her prescription function $\gamma_t^{l}$ to operate on her private information $x_{1:t}^{l}$ to produce her action $a_t^{l}$, i.e. $\gamma_t^{l}:(\cX^{l})^t\to \cP(\cA^{l})$ and $a_t^{l}\sim\gamma_t^{l}(\cdot|x_{1:t}^{l})$. And follower $j$ uses her prescription function $\gamma_t^{m,j}$ to operate on her private information $x_{1:t}^{m,j}$ to produce her action $a_t^{m,j}$, i.e. $\gamma_t^{m,j}:(\cX^{m,j})^t\to \cP(\cA^{m,j})$ and $a_t^{m,j}\sim\gamma_t^{m,j}(\cdot|x_{1:t}^{m,j})$. It is easy to see that for any $\sigma^{l,m}$ policy profile of the players, there exists an equivalent $\psi$ profile of the common agent (and vice versa) that generates the same control actions for every realization of the information of the players.

Here, we will consider Markovian common agent's policy as follows. We call a common agent's policy be of ``type $\theta$" if the common agent observes common belief $\underline{\pi_t}$, and generates prescription functions $\gamma_t := (\gamma_t^{l},\gamma_t^m) = \theta_t[\underline{\pi_t}]$. 
The Stackelberg leader uses her prescription function $\gamma_t^{l}$ to operate on her private information $x_t^{l}$ to produce her action $a_t^{l}$, i.e. $\gamma_t^{l}:\cX^{l}\to \cP(\cA^{l})$ and $a_t^{l}\sim\gamma_t^{l}(\cdot|x_{t}^{l})$. Similarly, the major follower $j$ uses her prescription function $\gamma_t^{m,j}$ to operate on her private information $x_t^{m,j}$ to produce her action $a_t^{m,j}$, i.e. $\gamma_t^{m,j}:\cX^{m,j}\to \cP(\cA^{m,j})$ and $a_t^{m,j}\sim\gamma_t^{m,j}(\cdot|x_{t}^{m,j})$.

Recall that we defined common marginal beliefs ${\underline{\pi_t}} = \{ \pi_t^{l}, \pi_t^{m,j}, \}_{ j=1\ldots M }$ , where $\pi_t^{l},\pi_t^{m,j}, j=1\ldots M$, where $\pi_t^{l}(x_t^{l}) = P^{\sigma^{l},\sigma^{m}}(x_t^{l}|a_{1:t-1}^{l,m})$,\\ $\pi_t^{m,j}(x_t^{m,j}) = P^{\sigma^{l},\sigma^{m}}(x_t^{m,j}|a_{1:t-1}^{l,m})$ i.e. $\pi_t^{l}$ is the common belief on the Stackelberg leader's state $x_t^{l}$ given the common information $(a_{1:t-1}^{l,m})$ and similarly $\pi_t^{m,j}$ is the common belief on the major follower $j$'s state $x_t^{m,j}$ given the common information $(a_{1:t-1}^{l,m})$. 
In the following lemma, we show that the belief $\underline{\pi_t}$ can be updated using Bayes' rule.
\begin{lemma}
There exists  functions $F^{l}, F^{m,j}$ for $ j=1\ldots M$, independent of the strategy $\theta$ such that 
\eq{
\pi^{l}_{t+1}&=F^{l}(\underline{\pi}_t^{l,m},\gamma^{l,m}_t,A^{l,m}_t)\\
\pi^{m,j}_{t+1}&=F^{m,j}(\underline{\pi}_t^{l,m},\gamma^{l,m}_t,A^{l,m}_t)
}
Combining the above two we also say
\eq{
\pi^{l,m}_{t+1}=\underline{F}(\underline{\pi}_t^{l,m},\gamma^{l,m}_t,a_t^{l,m})
}
\end{lemma}
\begin{proof}
We only consider the proof of the update of $\pi_t^{l}$ and the proof of the update of $\pi_t^{m,j}$ is similar which is skipped.
\eq{
\pi_{t+1}^{l}(x_{t+1}^{l}) &=  P^{\theta}(x_{1+t}^{l}|a_{1:t}^{l,m})
}
\eq{
&=\frac{\displaystyle \sum_{x_t^{l,m},a^{l,m}_t}\underline{\pi}_t^{l,m}(x^{l,m}_t)\gamma_t^{l,m}(a^{l,m}_t|x^{l,m}_t)Q^{l}(x_{t+1}^{l}|x^{l,m}_t,a^{l,m}_t)}{\displaystyle \sum_{x_t^{l,m}}\underline{\pi}_t^{l,m}(x_t^{l,m})\gamma_t^{l,m}(a^{l,m}_t|x^{l,m}_t)}\\
&=\frac{\displaystyle \sum_{x^{l,m}_t,a^{l,m}_t}\underline{\pi}_t(x^{l,m}_t)\gamma_t^{l,m}(a^{l,m}_t|x^{l,m}_t)Q^{l}(x_{t+1}^{l}|x^{l,m}_t,a^{l,m}_t)}{\displaystyle \sum_{x^{l,m}_t}\underline{\pi_t}(x^{l,m}_t)\gamma_t^{l,m}(a^{l,m}_t|x^{l,m}_t)}
\label{eq:piupdate}
}
\end{proof}


\section{Algorithm for SPE computation}
\label{sec:Result}
In the next section, we design an algorithm to compute SPE of the game.

\subsection{Backward Recursion}

In this section, we define an equilibrium generating function $\theta=(\theta^{l}_t,\theta_t^{m,j})_{i\in\{1\ldots N\}, j\in\{1\ldots M, \},t\in[T]}$, where $\theta^{l}_t :  (\mathcal{P}(\cX^{l}))\times(\prod_{j=1}^M\cP(\cX^{m,j}))\times \cP(\cX^{m,j}) \to \big\{\cX^{{l}} \to \mathcal{P}(\cA^{l}) \big\}$, and a sequence of functions $(V_t^{l}, V_t^{m,j})_{t\in \{ 1,2, \ldots T+1\}}$, where $V^{l}_t :  (\mathcal{P}(\cX^{l}))\times(\prod_{j=1}^M\cP(\cX^{m,j}))\times \cX^{l} \to \mathbb{R}$, $V^{m,j}_t :  (\mathcal{P}(\cX^{l}))\times(\prod_{j=1}^M\cP(\cX^{m,j}))\times  \cX^{m,j} \to \mathbb{R}$, in a backward recursive way, as follows.
 
\begin{itemize}
\item[1.] Initialize $\forall \underline{\pi}_{T+1}\in\mathcal{P}(\mathcal{X}^{l})\times\prod_{j=1}^M\mathcal{P}(\mathcal{X}^{m,j}), , j=1\ldots M, x_{T+1}^{l}\in \cX^{l}, x_{T+1}^{m,j}\in \cX^{m,j},$
\eq{
V^{l}_{T+1}(\underline{\pi}_{T+1},x_{T+1}^{l}) &\defeq 0.   \label{eq:VT+1}\\
V^{m,j}_{T+1}(\underline{\pi}_{T+1},x_{T+1}^{m,j}) &\defeq 0.   
}
\item[2.] For $t = T,T-1, \ldots 1, \ \forall \underline{\pi}_{t}\in\mathcal{P}(\mathcal{X}^{l})\times\prod_{j=1}^M\mathcal{P}(\mathcal{X}^{m,j})$,
$\gamma_t^{l},\gamma_t^{m,-j}\gamma_t^{m,j}$, define\\ $\bar{BR}_t^{m,j}(\underline{\pi_t},\gamma_t^{l},\gamma_t^{m,-j})$ as follows, 
\eq{
 &\bar{BR}_t^{m,j}(\underline{\pi_t},\gamma_t^{l},\gamma_t^{m,-j}) :=\big\{ \tgamma_t^{m,j}: \forall x_t^{m,j}\in \cX^{m,j}, \tgamma_t^{m,j}(\cdot|x_t^{m,j})\nn\\
 &\in  \arg\max_{\gamma^{m,j}_t(\cdot|x_t^{m,j})}\E^{\gamma^{m,j}_t(\cdot|x_t^{m,j}) {\gamma}^{l}_t,\gamma_t^{m,-j},\,\underline{\pi_t}} \nn\\
 & 
\big\{ R_t^{m,j}(X^{l,m}_t,A^{l,m}_t) +\delta V_{t+1}^{m,j}(\underline{F}(\underline{\pi_t},\gamma_t^{l},\tgamma_t^{m,j},\gamma_t^{m,-j},A^{l,m}_t), X^{m,j}_{t+1}) \big\lvert \underline{\pi_t}, x_t^{m,j} \big\}  \big\}, \label{eq:m_FP_mj}
}
where the expectation in (\ref{eq:m_FP_mj}) is with respect to random variables $(X^{l,m}_t,A^{l,m}_t)$ through the measure $\underline{\pi_t}(x^{l,m}_t) {\gamma}^{l,m}_t(a^{l,m}_t|x^{l,m}_t)$.
Then let for all $\underline{\pi_t}$, $\theta[\underline{\pi_t}] =(\tgamma_t^{l,m})$ is a solution of the following fixed-point equation (if it exists). For all $x_t^{l}$
\seq{

\eq{
&\tgamma_t^{l} \in \arg\max_{\gamma_t^{l}(\cdot|x_t^l)}\nn\\
&\E^{ {\gamma}^{l}_t(\cdot|x_t^l),\hat{\gamma}_t^m} \big\{ R_t^{l}(X^{l,m}_t,A^{l,m}_t) +\delta V_{t+1}^{l}(\underline{F}(\underline{\pi_t},\tgamma_t^{l},\hat{\gamma}^{m}_t,A^{l,m}_t),X_{t+1}^{l})|\underline{\pi_t},x^{l}_t\big\},  \label{eq:P2}\\
&\text{where } \hat{\gamma}_t^m \in \bar{BR}_t^{m}(\underline{\pi_t}, \tgamma_t^{l},\hat{\gamma}_t^{m})
}
where the above expectation is defined with respect to random variables \\$(X_t^{l,m},A^{l,m}_t)$ through the measure $\underline{\pi_t}(x^{l,m}_t) {\gamma}^{l}_t(a^{l}_t|x_t^{l})\hat{\gamma}^{m}_t(a^{m}_t|x_t^{m})\\Q^{l,m}(x_{t+1}^{l,m}|x_t^{l,m},a_t^{l,m})$. 
}
Let $(\tgamma_t^{l},\tgamma_t^{m})$ be a tuple of solution of the above operation. Then set $\forall j, x_{t}^{l},x_t^{m,j}$,
 \seq{
 \label{eq:Vdef}
  \eq{
   &V^{l}_{t}(\underline{\pi_t},x_t^{l}) \defeq  \nn\\
   &\;\E^{ \tilde{\gamma}^{l,m}_t}\big\{ {R}_t^{l} (X^{l,m}_T,A_t^{l,m}) + \delta V_{t+1}^{l} (\underline{F}(\underline{\pi_t},\tgamma^{l,m}_t,A^{l,m}_t),X_{t+1}^{l})\big\lvert \underline{\pi_t},x_t^{l} \big\}\\
   &V^{m,j}_{t}(\underline{\pi_t},x_t^{m,j}) \defeq \nn\\
   &\;\E^{ \tilde{\gamma}^{l,m}_t}\big\{ {R}_t^{m,j} (X^{l,m}_T,A_t^{l,m}) + \delta V_{t+1}^{m,j} (\underline{F}(\underline{\pi_t},\tgamma^{l,m}_t,A^{l,m}_t),X_{t+1}^{m,j})\big\lvert \underline{\pi_t},x_t^{m,j} \big\}
   }
   }
   \end{itemize}

Based on $\theta$ defined in the backward recursion above, we now construct a set of strategies $\tsigma$ through forward induction as follows. 

For $t =1,2 \ldots T,j, \underline{\pi_t},x_{1:t}^{l} \in(\cX^{l})^t,x_{1:t}^{m} \in(\cX^{m})^t,a_{1:t-1}^{l}\in(\cA^{l})^{t-1},a_{1:t-1}^{m}\in(\cA^{m})^{t-1}$
\eq{
\pi_1(x_1^{l,m}) &:= Q^{l,m}(x^{l}_1,x^{m}_1)\nn\\
\tsigma_{t}^{m,j}(a_{t}^{m,j}|a_{1:t-1}^{l,m}, x_{1:t}^{m,j}) &:= \theta_{t}^{m,j}[\underline{\pi_t}](a^{m,j}_{t}| x_{t}^{m,j})\\
\tsigma_{t}^{l}(a_{t}^{l}|a_{1:t-1}^{l},x_{1:t}^{l}) &:= \theta_{t}^{l}[\underline{\pi_t}](a^{l}_{t}|x_t^{l})  \\
\tsigma_{t}^{m,j}(a_{t}^{m,j}|a_{1:t-1}^{m,j},x_{1:t}^{m,j}) &:= \theta_{t}^{m,j}[\underline{\pi_t}](a^{m,j}_{t}|x_t^{m,j})  \\
\pi_{t+1} &= \underline{F}(\underline{\pi_t},\theta_t^{l,m}[\underline{\pi_t}],a^{l,m}_t)
}

\begin{theorem}
\label{Thm:Main}
A strategy profile $\tsigma$, as constructed through backward-forward recursion algorithm above is an SPE of the game
\end{theorem}
\begin{proof}
We will prove this theorem in two parts. 

In Part 1 for the major follower j, we prove that $\tsigma^{m,j} \in BR^{m,j}(\tsigma^{l},\tsigma^{m,-j},\tsigma_t^{m,j})$ i.e. $\ \forall \ t\in[T]$, $\forall \sigma^{m,j}, a_{1:t-1}^{l,m},x_{1:t}^{m,j}$
\eq{
 &\E^{\tsigma_{t:T}^{l,m},\underline{\pi_t}} \big\{ \sum_{n=t}^T \delta^{n-t}R_n^{m,j}(X^{l,m}_n,A^{l,m}_n) |\underline{\pi_t},a_{1:t-1}^{l,m},x_{1:t}^{m,j}\big\} \geq \nn\\
 &\E^{\tsigma_{t:T}^{l},\sigma_{t:T}^{m,j},\tsigma_{t:T}^{m,-j},\underline{\pi_t}} \big\{ \sum_{n=t}^T \delta^{n-t}R_n^{m,j}(X^{l,m}_n,A^{l,m}_n) |\underline{\pi_t},a_{1:t-1}^{l,m},x_{1:t}^{m,j}\big\} \label{eq:Thm_f}.
}

In Part 2 for the leader, we show that $\forall t,\sigma^{l},a_{1:t-1}^{l,m},x_{1:t}^{l}$
\eq{
&\E^{\tsigma_{t:T}^{l},{\tsigma}_{t:T}^{m},\underline{\pi_t}} \big\{ \sum_{n=t}^T \delta^{n-t}R_n^{l}(X_n^{l,m},A_n^{l,m}) |\underline{\pi_t},a_{1:t-1}^{l,m},x_{1:t}^{l}\big\} \geq\nn\\
 &\E^{\sigma_{t:T}^{l},\hat{\sigma}_{t:T}^{m},\hat{\sigma}_{t:T}^{m,j},\underline{\pi_t}} \big\{ \sum_{n=t}^T \delta^{n-t}R_n^{l}(X_n^{l,m},A_n^{l,m}) |\underline{\pi_t},a_{1:t-1}^{l,m},x_{1:t}^{l}\big\},\\
 &\text{where } \hat{\sigma}^{m,j}\in BR^{m,j}(\sigma_{t:T}^{l},\hat{\sigma}_{t:T}^{m}),\hat{\sigma}^m\in BR^m(\sigma_{t:T}^{l},\hat{\sigma}_{t:T}^{m},\hat{\sigma}_{t:T}^{m,j})
}
where $\tsigma^{m,j}\in BR^{m,j}(\tsigma^{l},\tsigma^{m})$, as shown in Part 1.

Combining the above parts prove the above result. The proof is presented in Appendix~C.
\end{proof}

In the following, we show that every Stackelberg mean field equilibrium can be found using the above backward recursion. This also enables us to comment on the existence of the solution of the fixed-point equation~\eqref{eq:P2}.

\begin{theorem}
Suppose there exists an SPE $(\tilde{\sigma}^{l},\tsigma^m)$ that is a solution of the fixed point equation defined in Definition~\ref{def:ME}. Then there exists an equilibrium generating function $\theta$ that satisfies \eqref{eq:P2} in backward recursion $\forall \underline{\pi_t}$
	such that  $(\tilde{\sigma}^{l},\tsigma^m)$ is defined through forward recursion using $\theta$. This also implies that there exists a solution of~\eqref{eq:P2} for each time $t$.
\end{theorem}
\begin{proof}
Suppose there exists an SPE $(\tilde{\sigma}^{l,m})$ of the game. The proof in Appendix~\ref{app:Proof_Exist} show that all SPE can be found using backward/forward recursion. This proves that there exists a solution of~\eqref{eq:P2} for every $t$.
\end{proof}

\textbf{Remark:}
When leader is social welfare maximizing, her utility can be given by 
\eq{
R^{l}(x_t^{l,m},a_t^{l,m}) &= \sum_{x_t^{l,m},a_t^{l,m}}R^{m,j}(x^{l,m}_t,a^{l,m}_t)
}

\section{Special case 2: Infinite horizon case}
\label{sec:Infinite}
In this section we consider the case with infinite horizon. For this section we assume that the instantaneous rewards of the players $R^{l},R^{m,j}$ are absolutely bounded and do not depend on time.

We design an algorithm to compute SPE of the infinite horizon game as follows.

\subsection{Backward Recursion}

In this section, we define an equilibrium generating function $\theta=(\theta^{l},\theta^{m,j})_{ j\in\{1\ldots M, \}}$, where $\theta^{l} :  (\mathcal{P}(\cX^{l}))\times(\prod_{j=1}^M\cP(\cX^{m,j})) \to \big\{\cX^{{l}} \to \mathcal{P}(\cA^{l}) \big\}$, and a sequence of functions $(V^{l}, V^{m,j})$, where $V^{l} :  (\mathcal{P}(\cX^{l}))\times(\prod_{j=1}^M\cP(\cX^{m,j}))\times \cX^{l} \to \mathbb{R}$, $V^{m,j} :  (\mathcal{P}(\cX^{l}))\times(\prod_{j=1}^M\cP(\cX^{m,j}))\times  \cX^{m,j} \to \mathbb{R}$, in a backward recursive way, as follows. 
\begin{itemize}
\item[1.] 

\seq{
Define $\forall \underline{\pi}_{t}\in\mathcal{P}(\mathcal{X}^{l})\times\prod_{j=1}^M\mathcal{P}(\mathcal{X}^{m,j})$, $\gamma_t^{l},\gamma_t^{m,-j}$, define $\bar{BR}_t^{m,j}(\underline{\pi},\gamma^{l},\gamma^{m,-j})$ as follows, 
\eq{
 &\bar{BR}^{m,j}(\underline{\pi},\gamma^{l},\gamma^{m,-j}) :=\big\{ \tgamma^{m,j}: \forall x^{m,j}\in \cX^{m,j}, \tgamma^{m,j}(\cdot|x^{m,j})\nn\\
 &\in  \arg\max_{\gamma^{m,j}(\cdot|x^{m,j})}\E^{\gamma^{m,j}(\cdot|x^{m,j}) {\gamma}^{l}\gamma_t^{m,-j},\underline{\pi}} \nn\\
 & 
\big\{ R^{m,j}(X^{l,m},A^{l,m}) +\delta V^{m,j}(\underline{F}(\underline{\pi},\gamma^{l},\tgamma^{m,j},\gamma_t^{m,-j},A^{l,m}), X^{m,j,'}) \big\lvert \underline{\pi},x^{m,j} \big\}  \big\}, \label{eq:m_P_mj3}
}
where the expectation in (\ref{eq:m_P_mj3}) is with respect to random variables $(X^{l,m},A^{l,m})$ through the measure $\underline{\pi}(x^{l,m}) {\gamma}^{l,m}(a^{l,m}|x^{l,m})$.
}
Then let for all $\underline{\pi}$, $\theta[\underline{\pi}] =(\tgamma^{l,m})$ is a solution of the following fixed-point equation (if it exists),
\seq{
\eq{
\tgamma^{m} &\in \bar{BR}^{m}(\underline{\pi},\tgamma^{l},\tgamma^{m}) \label{eq:P1b}\\
\tgamma^{l} &\in \arg\max_{\gamma^{l}(\cdot|x^l)}\E^{ {\gamma}^{l},\hat{\gamma}^m} \big\{ R^{l}(X^{l,m},A^{l,m}) +\delta V^{l}(\underline{F}(\underline{\pi},\tgamma^{l},\hat{\gamma}^{m},A^{l,m}),X^{l,'})|\underline{\pi},x^{l}\big\},  \label{eq:P2b}\\
&\text{where }  \hat{\gamma}^m \in \bar{BR}^{m}(\underline{\pi},\gamma^{l},\hat{\gamma}^{m}),
}
where the above expectation is defined with respect to random variables \\$(X^{l,m},A^{l,m}_t)$ through the measure $\underline{\pi}(x^{l,m}) {\gamma}^{l}(a^{l}|x^{l})\hat{\gamma}^{m}(a^{m}|x_t^{m})$\\$Q^{l,m}(x^{l,m,'}|x^{l,m},a^{l,m})$. 
}
And $\forall j,x^{l},x^{m,j}$,
 \seq{
 \label{eq:Vdef}
  \eq{
   &V^{l}(\underline{\pi},x^{l}) \defeq  \;\E^{\tilde{\gamma}^{l,m}}\big\{ {R}^{l} (X^{l,m},A^{l,m}) + \delta V^{l} (\underline{F}(\underline{\pi},\tgamma^{l,m},A^{l,m}),X^{l,'})\big\lvert \underline{\pi},x^{l} \big\}\\
   &V^{m,j}(\underline{\pi},x^{m,j}) \defeq  \nn\\
   &\;\E^{ \tilde{\gamma}^{l,m}}\big\{ {R}^{m,j} (X^{l,m},A^{l,m}) + \delta V^{m,j} (\underline{F}(\underline{\pi},\tgamma^{l,m},A^{l,m}),X^{m,j,'})\big\lvert \underline{\pi},x^{m,j} \big\}
   }
   }
   \end{itemize}

Based on $\theta$ defined in the backward recursion above, we now construct a set of strategies $\tsigma$ through forward induction as follows. 

For $t =1,2 \ldots \infty,j, \underline{\pi_t},x_{1:t}^{l} \in(\cX^{l})^t,x_{1:t}^{m} \in(\cX^{m})^t,a_{1:t-1}^{l}\in(\cA^{l})^{t-1},a_{1:t-1}^{m}\in(\cA^{m})^{t-1}$
\eq{
\pi_1(x_1^{l,m}) &:= Q^{l,m}(x^{l}_1,x^{m}_1)\nn\\
\tsigma_{t}^{l}(a_{t}^{l}|a_{1:t-1}^{l},x_{1:t}^{l}) &:= \theta_{t}^{l}[\underline{\pi_t}](a^{l}_{t}|x_t^{l})  \\
\tsigma_{t}^{m,j}(a_{t}^{m,j}|a_{1:t-1}^{l,m},x_{1:t}^{m,j}) &:= \theta_{t}^{m,j}[\underline{\pi_t}](a^{m,j}_{t}|x_t^{m,j})  \\
\pi_{t+1} &= \underline{F}(\underline{\pi_t},\theta_t^{l,m}[\underline{\pi_t}],a^{l,m}_t)
}

\begin{theorem}
\label{Thm:Main}
A strategy profile $\tsigma$, as constructed through backward-forward recursion algorithm above is an SPE of the game
\end{theorem}
\begin{proof}
The proof is similar to the extension of finite horizon problems to infinite horizon problems in standard stochastic control problems and for now we omit the proof.
\end{proof}
\section{Examples}
\subsection{Security Example}
\label{sec:Example}
\begin{figure}[htbp] 
   \centering
   \begin{subfigure}{.5\textwidth}
  \centering
  \includegraphics[width=2in]{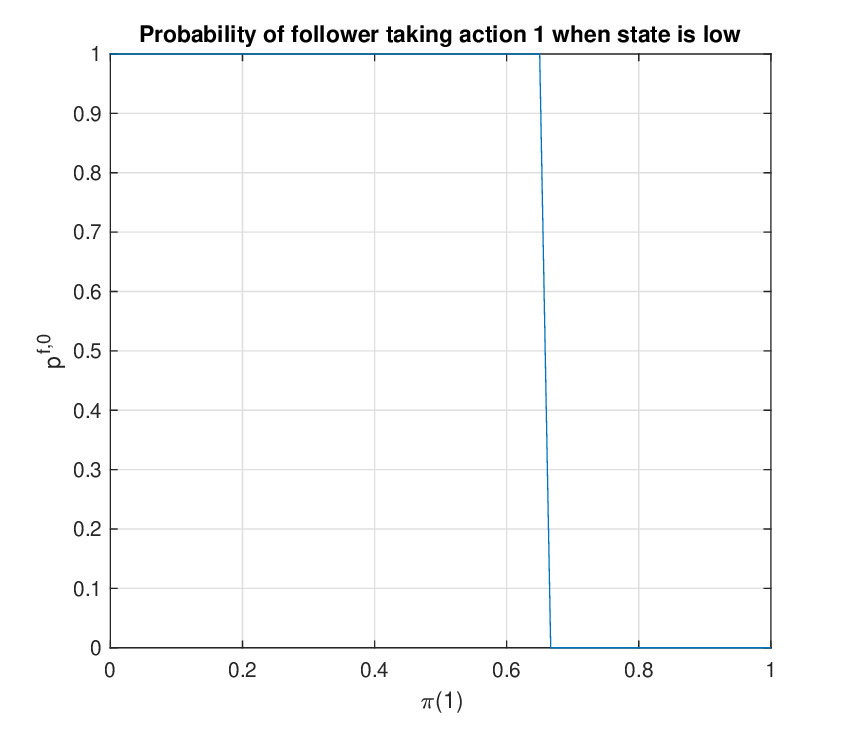} 
 \caption{Low}
  \label{fig:sub1}
\end{subfigure}%
\begin{subfigure}{.5\textwidth}
  \centering
 \includegraphics[width=2in]{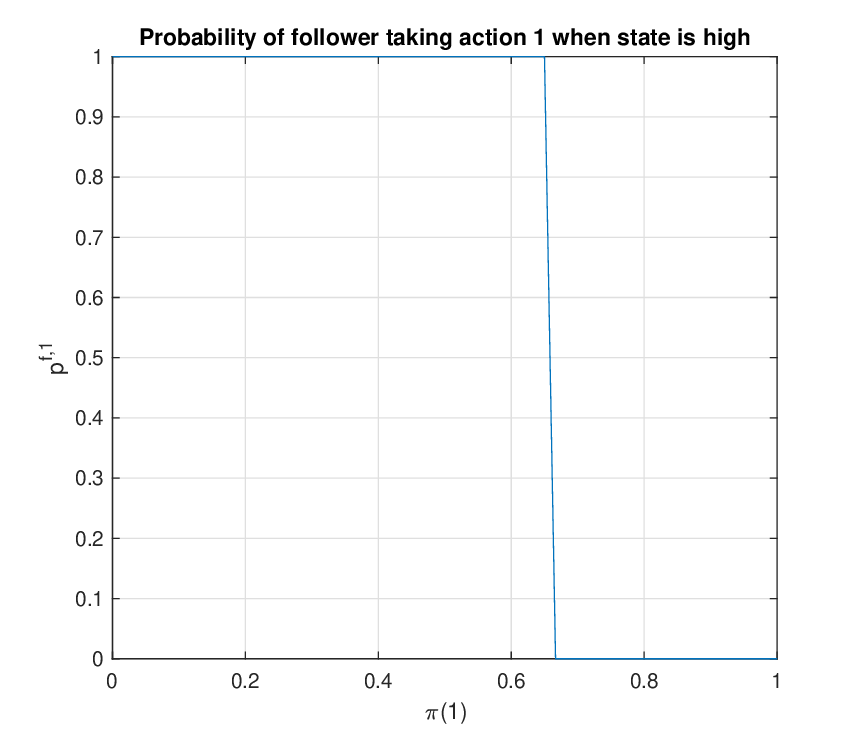} 
  \caption{High}
  \label{fig:sub2}
\end{subfigure}
     \caption{Probability of follower taking action 1 when its state is low and high}
\label{fig:example}
\end{figure}

\begin{figure}[htbp] 
   \centering
   \begin{subfigure}{.5\textwidth}
  \centering
  \includegraphics[width=2in]{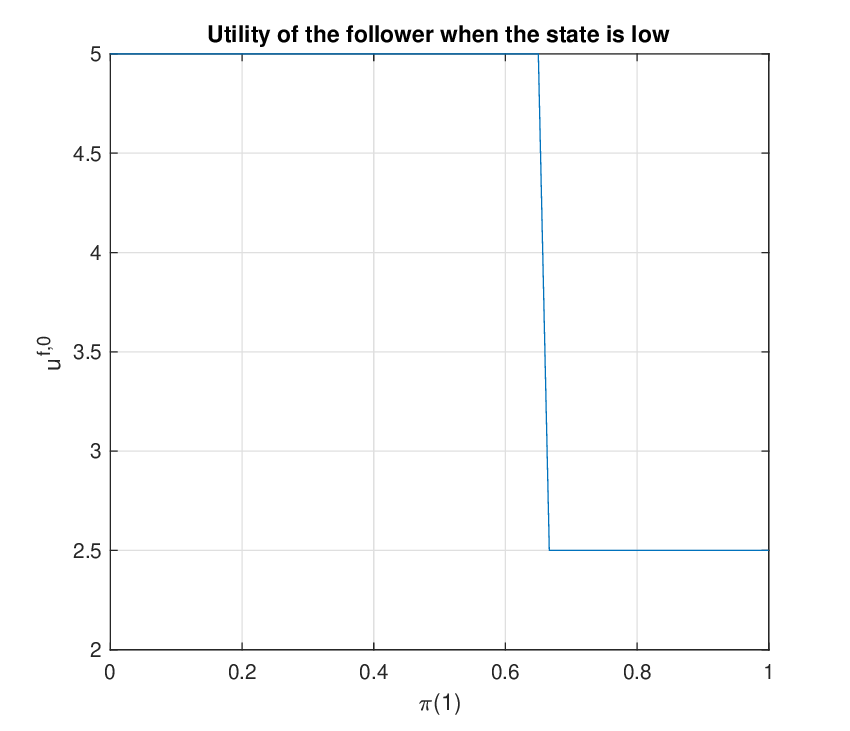} 
 \caption{Low}
  \label{fig:sub1}
\end{subfigure}%
\begin{subfigure}{.5\textwidth}
  \centering
 \includegraphics[width=2in]{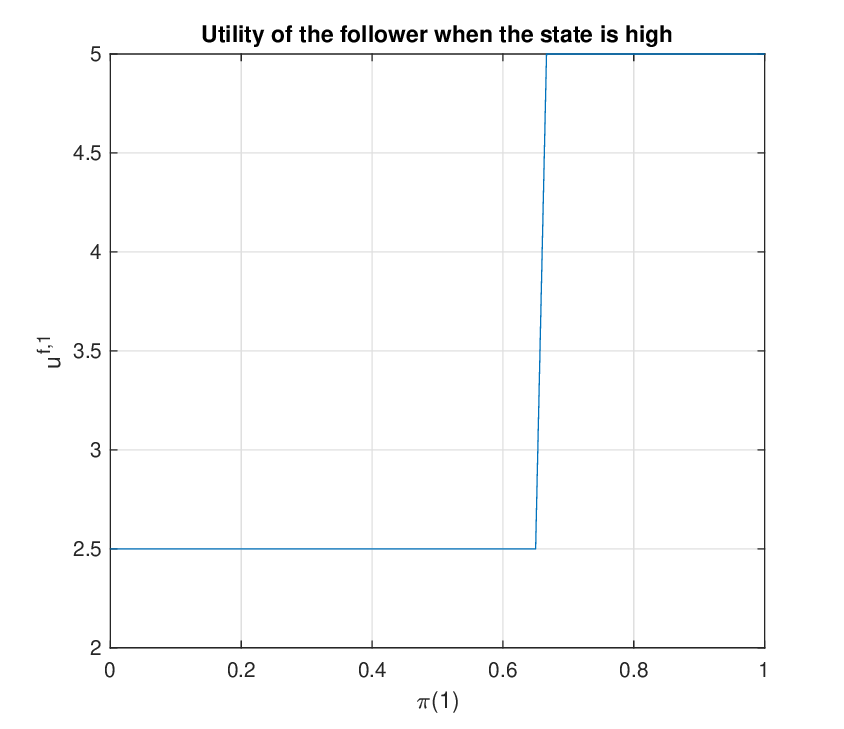} 
  \caption{High}
  \label{fig:sub2}
\end{subfigure}
     \caption{Utility of the follower when its state is low and high}
\label{fig:example2}
\end{figure}

\begin{figure}[htbp] 
   \centering
   \begin{subfigure}{.5\textwidth}
  \centering
  \includegraphics[width=2in]{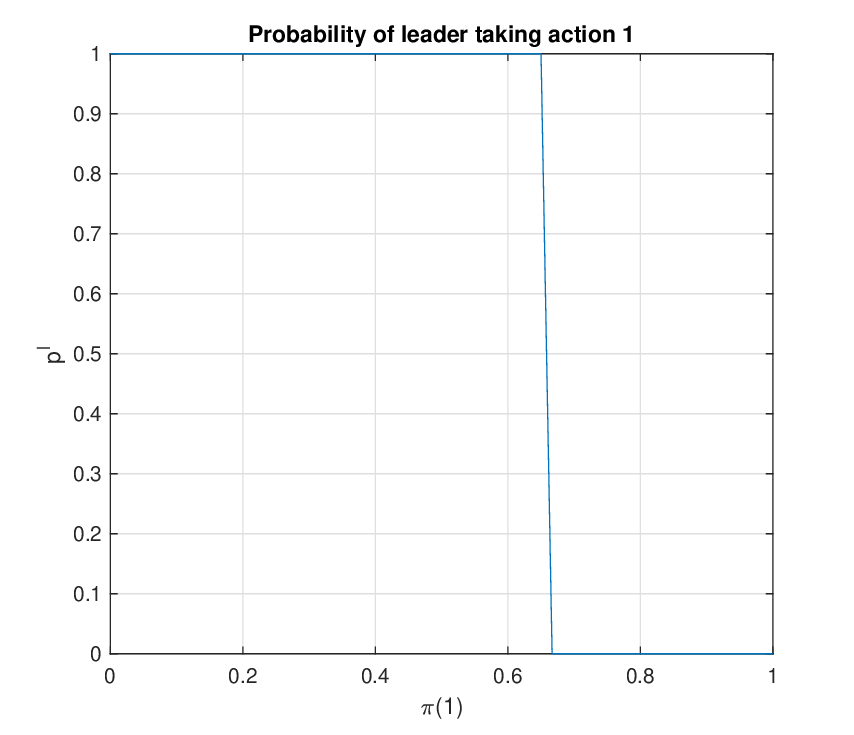} 
  \caption{Probability of leader taking action 1 }
  \label{fig:sub1}
\end{subfigure}%
\begin{subfigure}{.5\textwidth}
  \centering
 \includegraphics[width=2in]{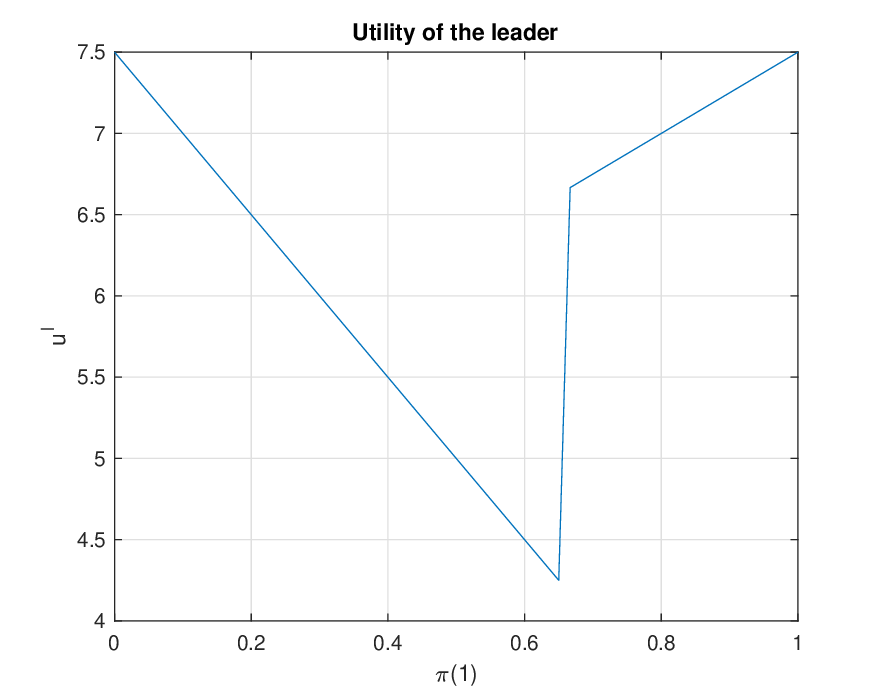} 
   \caption{Utility of leader }
     \label{fig:sub2}
\end{subfigure}
\caption{}
  \label{fig:example4}
\end{figure}

In this section, we consider a repeated Stackelberg game as a security example. We assume that $\cX^{m,j} = \cA^l = \cA^{m,j}=\{0,1\}, \cX^l = \phi$ and type of the leader is static i.e. $Q(x_{t+1}|x_t,a_t) = \mathbb{1}(x_{t+1}=x_t)$. We assume $\delta = 0.6$.
Let $p^l = \gamma^l(1)$, $p^{f,0} =\gamma^{m,j}(1|0) $ and  $p^{f,1} =\gamma^{m,j}(1|1)$ and the rewards of the players are given in Table~I below.
\begin{center}

\begin{table}
\caption{Game matrix for X=0 \& X=1}
\quad \quad\quad \quad
\begin{tabular}{|cc|c|c|}
\hline
      X=0        &                 & follower    & follower\\
              &                 & $A1$               & $A2$ \\
\hline
leader  & $D1$               & $(2,1)$           & $(4,0)$ \\
\hline          
   leader           & $D2$           & $(1,0)$         & $(3,2)$ \\
\hline
\end{tabular}
\begin{tabular}{|cc|c|c|}
\hline
      X=1        &                 & follower    & follower\\
              &                 & $A1$               & $A2$ \\
\hline
leader  & $D1$               & $(3,2)$           & $(2,0)$ \\
\hline          
   leader           & $D2$           & $(0,1)$         & $(1,1)$ \\
\hline
\end{tabular}
\end{table}
\end{center}
The equilibrium strategies and value functions are provided in igures~\ref{fig:example}--\ref{fig:example4}. Interestingly, the equilibrium strategies of the players are pure strategies that exhibit ``complementary discontinuities"~\cite{DaMa86i,DaMa86ii}.

\subsection{Beeps}
We consider the discrete time version of Example in~\cite{El17}.
Let $x_t\in\{0,1\}$ be a Markov process privately observed by a leader who sends a signal at each time $t$ to a follower, who upon receiving this signal maintains a belief $\pi_t$ on the state $x_t$. $x_t$ evolves as
\eq{
x_{t+1} = x_t + (1-x_t)m_t
}
where $m_t = 1$ with probability p independent across time. The process starts from 0 and 1 is an absorbing state. At each time $t$ the leader sends a signal $s_t$ to the follower about $x_t$ such that $s_t = \sigma^l_t(x_{1:t})$. Let $\pi_t$ be a common belief maintained by both the leader and the follower where for all $x_t$
\eq{
\pi_{t}(x_t) = P^{\sigma^l}(x_t|s_{1:t})
}
Thus
\eq{
\pi_{t+1}(x_{t+1}) = \frac{\sum_{x_t}P(x_t,x_{t+1},s_{t+1}|s_{1:t})}{\sum_{x_t,x_{t+1}}P(x_t,x_{t+1},s_{t+1}|s_{1:t})}\\
\frac{\sum_{x_t}\pi_t(x_t)Q(x_{t+1}|x_t)\gamma_{t+1}^l(s_{t+1}|x_{t+1})}{\sum_{x_t,x_{t+1}}\pi_t(x_t)Q(x_{t+1}|x_t)\gamma_{t+1}^l(s_{t+1}|x_{t+1})}\\
}
\eq{
\pi_{t+1} = F(\pi_t,\gamma^l_{t+1},s_{t+1})
}
The leader's payoff is
\eq{
R^l(\pi_t) = \lb{1 \text{ if } \pi_t\leq p^* \\0 \text{ if } \pi_t> p^*
}
}
Since the follower does not influence the state, she is myopic and her best response is built in the above payoff of the leader. The leader wants to find the policy that it can commit to such that it maximizes her infinite horizon cumulative discounted payoff

From Section~\ref{sec:Infinite}, one needs to solve 

\eq{
\tgamma_t^{l} &\in \arg\max_{\gamma^{l}(\cdot|x^l)}\E^{ {\gamma}^{l}} \big\{ R^{l}(\pi) +\delta V^{l}(\underline{F}(\underline{\pi},\tgamma^{l},S),X^{l,'})|\underline{\pi},x^{l}\big\}\\
V^{l}(\underline{\pi},x^{l}) \defeq  &\;\E^{\tilde{\gamma}^{l}}\big\{ R^{l}(\pi) + \delta V^{l} (\underline{F}(\underline{\pi},\tgamma^{l},S),X^{l,'})\big\lvert \underline{\pi},x^{l} \big\}
}
In the following, we show the plot of the utility of the leader in Figure 4. Note that it is different from the utility of the leader obtained in~\cite{El17} since they average out the private state of the player and poses the problem in common belief, which can then be posed as an optimization problem. However, in our case, we condition on the private state and thus rather get a fixed-point equation shown above.

\begin{figure}[htbp] 
   \centering
  \includegraphics[width=3in]{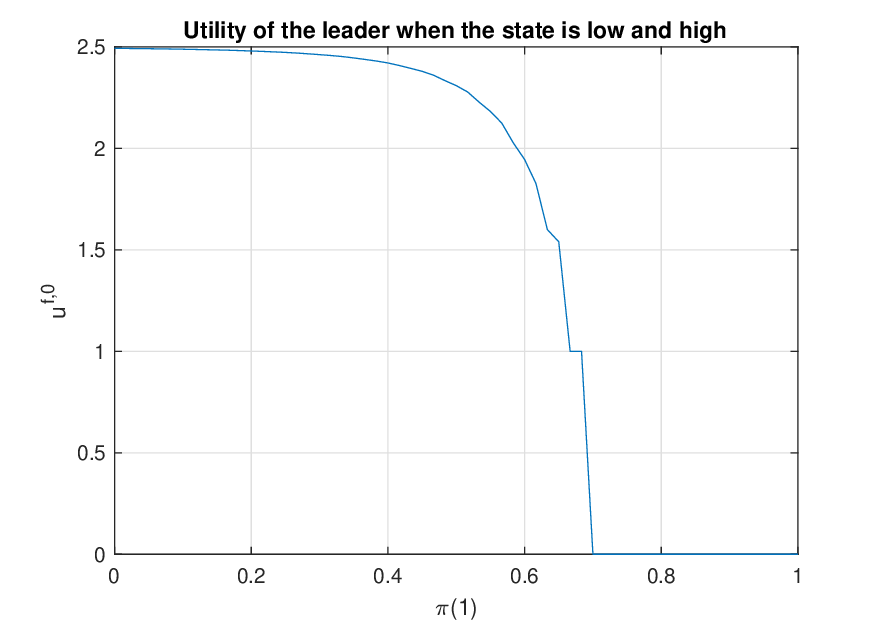} 
  \caption{Utility of the leader for both cases when the state is low and high (Example beeps)}
\end{figure}


\section{Conclusion}
\label{sec:Conclusion}
In this paper, we study a general stochastic Stackelberg game with single leader where both the followers and the leader have private types that evolves as conditionally independent controlled Markov process, conditioned on action history. We present a novel dynamic programming like methodology to sequentially decompose the problem of computing stochastic Stackelberg equilibrium for these games. Based on this algorithm we study a repeated security game where we numerically compute the equilibrium policies. In general, this algorithm can further increase the applicability of Stackelberg security games in dynamic security settings and in dynamic mechanism design where a leader commits to a policy and the follower best responds to it.

\section{}
\label{app:0}
\begin{claim}
	For any policy profile $g$ and $\forall t$,
	\eq{
	\p^{\sigma}(x_{1:t}^{l,m}|a_{1:t-1}) =  \p^{\sigma^{l}}(x_{1:t}^{l}|a_{1:t-1})\p^{\sigma^m}(x_{1:t}^{m}|a_{1:t-1})
	}
	\label{claim:CondInd}
	\end{claim}
	\begin{proof}
	\seq{
	\eq{
	&\p^{\sigma}(x_{1:t}|a_{1:t-1})= \frac{\p^{\sigma}(x_{1:t},a_{1:t-1})}{\sum_{\bar{x}_{1:t}} \p^{\sigma}(\bar{x}_{1:t},a_{1:t-1})}
	}
	Here, we will take numerator and the denominator separately.
	\eq{
	&Nr =  \left(Q_1^{l}(x^{l}_1)\sigma^{l}_1(a_1^{l}|x_{1}^{l})\prod_{n=2}^t Q_n^{l}(x^{l}_{n}|a_{n-1}^{l,m},x^{l}_{n-1}) \sigma^{l}_n(a_n^{l}|a^{l,m}_{1:n-1},x_{1:n}^{l}) \right)\\
	&\times \left(Q_1^{m}(x^{m}_1)\sigma^{m}_1(a_1^{m}|x_{1}^{m})\prod_{n=2}^t Q_n^{m}(x^{m}_{n}|a_{n-1}^{l,m},x^{m}_{n-1}) \sigma^{m}_n(a_n^{m}|a^{l,m}_{1:n-1},x_{1:n}^{m}) \right)
	}
	and 
	\eq{
	Dr&=\sum_{x_{1:t}^l} =  \left(Q_1^{l}(x^{l}_1)\sigma^{l}_1(a_1^{l}|x_{1}^{l})\prod_{n=2}^t Q_n^{l}(x^{l}_{n}|a_{n-1}^{l,m},x^{l}_{n-1})) \sigma^{l}_n(a_n^{l}|a^{l,m}_{1:n-1},x_{1:n}^{l}) \right)\\
	&\times\sum_{x_{1:t}^m} \left(Q_1^{m}(x^{m}_1)\sigma^{m}_1(a_1^{m}|x_{1}^{m})\prod_{n=2}^t Q_n^{m}(x^{m}_{n}|a_{n-1}^{l,m},x^{m}_{n-1}) \sigma^{m}_n(a_n^{m}|a^{l,m}_{1:n-1},x_{1:n}^{m}) \right)
	}
	
	Thus
	\eq{
	\p^{\sigma}(x_{1:t}^{l},x_{1:t}^{m}|a_{1:t-1}) =  \p^{\sigma^{l}}(x_{1:t}^{l}|a_{1:t-1})\p^{\sigma^{m}}(x_{1:t}^{m}|a_{1:t-1})
	}
	}
	\end{proof}
	\section{}
	\label{app:gsm}
	For any player $i$ (leader or follower), we use the notation $g$ to denote a general policy of the form $A_t^i\sim g_t^i(\cdot|a_{1:t-1}, x_{1:t}^{i})$, notation $s$ to denote a policy of the form $A_t^i \sim s_t^i(\cdot|a_{1:t-1},x_t^i)$, and notation $m$ to denote a policy of the form $A_t^i \sim m_t^i(\cdot|\pi_t,x_t^i)$. It should be noted that since $\pi_t$ is a function of random variables $a_{1:t-1}$, $m$ policy is a special type of $s$ policy, which in turn is a special type of $g$ policy.

Using the agent-by-agent approach~\cite{Ho80}, we show in Lemma~\ref{fact:G2S} that any expected reward profile of the players that can be achieved by any general strategy profile $g$ can also be achieved by a strategy profile $s$.
\begin{lemma}
Given a fixed strategy $g^{-i}$ of all players other than player $i$ and for any strategy $g^i$ of player $i$, there exists a strategy $s^i$ of player $i$ such that $\forall t \in \mathcal{T}, x_t\in \cX, a_t\in \cA,$
\eq{P^{s^i g^{-i}}(x_t, a_t) &= P^{g^ig^{-i}}(x_t, a_t) \;\;\;\;\;
}
which implies $ J^{i,s^ig^{-i}} = J^{i,g^ig^{-i}}$.\label{fact:G2S}
\end{lemma}
\begin{proof}
	The proof is on the similar lines as the proof of Lemma~1 in~\cite{VaSiAn19} 
\end{proof}
Since any $s^i$ policy is also a $g^i$ type policy, the above lemma can be iterated over all players which implies that for any $g$ policy profile there exists an $s$ policy profile that achieves the same reward profile i.e., $(J^{i,s})_{i\in \cN} = (J^{i,g})_{i\in \cN}$.
In the following lemma, we show that the space of profiles of type $s$ is outcome-equivalent to the space of profiles of type $m$.
\begin{lemma}
\label{fact:L1}
For any given strategy profile $s$ of all players, there exists a strategy profile $m$ such that
\eq{P^m(x_t, a_t) &= P^s(x_t, a_t) \;\;\;\;\forall t \in \mathcal{T}, x_t\in \cX, a_t\in \cA ,
}
which implies $ (J^{i,m})_{i\in \cN} = (J^{i,s})_{i\in \cN} $.
\label{fact:S2M}
\end{lemma}
\begin{proof}
	The proof is on the similar lines as the proof of Lemma~2 in~\cite{VaSiAn19}
\end{proof}

\section{Part 1: Followers}
\label{b_app:P1}
\label{b_app:B}
\label{b_app:A}
\begin{proof}
We prove Theorem~\ref{Thm:Main} using induction and the results in \\
Lemma~\ref{b_lemma:2}, and \ref{b_lemma:1} proved in Appendix~\ref{b_app:B}. Let $\tsigma$ be the strategies computed by the methodology in Section~III.
\seq{
For the base case at $t=T$, $ a_{1:T-1}^{l,m},x_{1:T}^{m,j},\sigma^{m,j}$
\eq{
&\E^{\tsigma_T^{l},\tsigma_{T}^{m,j},\tsigma_T^{m,-j},\underline{\pi_t}}\big\{  R^{m,j}_T(X^{l,m}_T,A^{l,m}_T) \big\lvert \underline{\pi_t},a_{1:T-1}^{l,m},x_{1:T}^{m,j} \big\}
=
V^{m,j}_T(\underline{\pi_T}, x_T^{m,j})  \label{b_eq:T2a}\\
&\geq \E^{\tsigma_T^{l},\sigma_{T}^{m,j},\tsigma_T^{m,-j},\underline{\pi_T}} \big\{ R^{m,j}_T(X^{l,m}_T,A^{l,m}_T) \big\lvert\underline{\pi_t},  a_{1:T-1}^{l,m},x_{1:T}^{m,j} \big\},  \label{b_eq:T2}
}
}
where \eqref{b_eq:T2a} follows from Lemma~\ref{b_lemma:1} and \eqref{b_eq:T2} follows from Lemma~\ref{b_lemma:2} in Appendix~\ref{b_app:B}.

Let the induction hypothesis be that for $t+1$. Then $\forall t,  a_{1:t}^{l,m},x_{1:t+1}^{m,j} \in (\cX^{m,j})^{t+1}, \sigma^{m,j}$,
\seq{
\eq{
 \E^{\tsigma_{t+1:T}^{l,m},\underline{\pi}_{t+1} } \big\{ \sum_{n=t+1}^T \delta^{n-t-1}R^{m,j}_n(X^{l,m}_n,A^{l,m}_n) \big\lvert \pi_{t+1}, a_{1:t}^{l,m}, x_{1:t+1}^{m,j} \big\} \\
 \geq
  \E^{\tsigma_{t+1:T}^{l},\sigma_{t+1:T}^{m,j},\tsigma_{t+1:T}^{m,-j},\underline{\pi}_{t+1} } \big\{ \sum_{n=t+1}^T \delta^{n-t-1} R^{m,j}_n(X^{l,m}_n,A^{l,m}_n) \big\lvert \pi_
  {t+1}, a_{1:t}^{l,m}, x_{1:t+1}^{m,j} \big\}. \label{b_eq:PropIndHyp}
}
}
\seq{
Then $\forall a_{1:t-1}^{l,m}, x_{1:t}^{m,j}, \sigma^{m,j}$, we have
\eq{
&\E^{\tsigma_{t:T}^{l,m},\underline{\pi_t} } \big\{ \sum_{n=t}^T \delta^{n-t-1}R^{m,j}_n(X^{l,m}_n,A^{l,m}_n) \big\lvert\underline{\pi_t}, a_{1:t-1}^{l,m},x_{1:t}^{m,j} \big\} \nonumber \\
&= V^{m,j}_t(\underline{\pi_t}, x_t^{m,j})\label{b_eq:T1}\\
&\geq \E^{\tsigma_t^{l},\sigma_t^{m,j},\tsigma_t^{m,-j},\underline{\pi_t}} \big\{ R^{m,j}_t(X^{l,m}_t,A^{l,m}_t) +\nn\\
&\delta V^{m,j}_{t+1} (\underline{F}(\underline{\pi_t},\tgamma^{l,m}_t,A^{l,m}_t)),X_{t+1}^{m,j}) \big\lvert\underline{\pi_t}, a_{1:t-1}^{l,m}, x_{1:t}^{m,j} \big\}  \label{b_eq:T3}\\
&= \E^{\tsigma_t^{l},\sigma_t^{m,j},\tsigma_t^{m,-j},\underline{\pi_t} } \big\{ R^{m,j}_t(X^{l,m}_t,A^{l,m}_t) + \delta \E^{\tsigma_{t+1:T}^{l,m},\underline{F}(\underline{\pi_t},\tgamma_t^{l,m},A^{l,m}_t)} \nn\\
&\big\{ \sum_{n=t+1}^T \delta^{n-t-1}R^{m,j}_n(X^{l,m}_n,A^{l,m}_n) \big\lvert \underline{F}(\underline{\pi_t},\tgamma^{l,m}_t,A^{l,m}_t), x_{1:t}^{m,j},X_{t+1}^{m,j} \big\}  \big\vert\underline{\pi_t}, a_{1:t-1}^{l,m}, x_{1:t}^{m,j} \big\}  \label{b_eq:T3b}\\
&\geq \E^{\tsigma_t^{l},\sigma_t^{m,j},\tsigma_t^{m,-j},\underline{\pi_t} } \big\{ R^{m,j}_t(X^{l,m}_t,A^{l,m}_t) + \delta\E^{\tsigma_{t+1:T}^{l,m},\underline{F}(\underline{\pi_t},\tgamma^{l,m}_t,A_t^{l,m}) } \nn\\
&\big\{ \sum_{n=t+1}^T \delta^{n-t-1}R^{m,j}_n(X^{l,m}_n,A^{l,m}_n) \big\lvert \underline{F}(\underline{\pi_t},\tgamma^{l,m}_t,A^{l,m}_t), x_{1:t}^{m,j},X_{t+1}^{m,j}\big\} \big\vert\underline{\pi_t}, a_{1:t-1}^{l,m}, x_{1:t}^{m,j} \big\}  \label{b_eq:T4} \\
&= \E^{\tsigma_t^{l},\tsigma_t^{m,j},\sigma_t^{m,-j},\underline{\pi_t} } \big\{ R^{m,j}_t(X^{l,m}_t,A^{l,m}_t) + \delta\E^{\tsigma_{t:T}^{l,m},\underline{\pi_t} }   \nn\\
& \big\{ \sum_{n=t+1}^T \delta^{n-t-1}R^{m,j}_n(X^{l,m}_n,A^{l,m}_n) \big\lvert \underline{F}(\underline{\pi_t},\tgamma^{l,m}_t,A_t^{l,m}), x_{1:t}^{m,j},X_{t+1}^{m,j}\big\} \big\vert\underline{\pi_t}, a_{1:t-1}^{l,m}, x_{1:t}^{m,j} \big\}
\label{b_eq:T5}\\
&=\E^{\tsigma_{t:T}^{l,m},\underline{\pi_t} } \big\{ \sum_{n=t}^T \delta^{n-t}R^{m,j}_n(X^{l,m}_n,A^{l,m}_n) \big\lvert\underline{\pi_t}, a_{1:t-1}^{l,m},x_{1:t}^{m,j} \big\}  \label{b_eq:T6},
}
}
where \eqref{b_eq:T1} follows from Lemma~\ref{b_lemma:1}, \eqref{b_eq:T3} follows from Lemma~\ref{b_lemma:2}, \eqref{b_eq:T3b} follows from Lemma~\ref{b_lemma:1}, \eqref{b_eq:T4} follows from induction hypothesis in \eqref{b_eq:PropIndHyp} and \eqref{b_eq:T5} follows from Lemma~\ref{b_lemma:3}.
\end{proof}

\section{}
\label{b_app:B}
\label{b_app:lemmas}
\begin{lemma}
\label{b_lemma:2}
Let $\tsigma$ be the strategies computed by the methodology in Section~III. 
Then $\forall t\in [T],a_{1:t-1}^{l,m}, x_{1:t}^{m,j}, \sigma^{m,j}_t$
\eq{
&V^{m,j}_t(\underline{\pi_t},x_t^{m,j})\geq \nn\\
& \E^{\tsigma_t^{l},\sigma_t^{m,j},\tsigma_t^{m,-j},\underline{\pi_t}} \big\{ R^{m,j}_t(X^{l,m}_t,A^{l,m}_t) + \delta V^{m,j}_{t+1} (\underline{F}(\underline{\pi_t},\tgamma^{l,m}_t,A^{l,m}_t), X_{t+1}^{m,j}) \big\lvert \underline{\pi_t}, a_{1:t-1}^{l,m}, x_{1:t}^{m,j} \big\}\label{b_eq:lemma2}
}
\end{lemma}

\begin{proof}
We prove this lemma by contradiction. Suppose the claim is not true for $t$. This implies $\exists  \hat{\sigma}_t^{m,j},\hat{a}_{1:t-1}^{l,m}$, $ \hat{x}_{1:t}^{m,j}$ such that
\eq{
&\E^{\tsigma_t^{l},\hat{\sigma}_t^{m,j},\tsigma_t^{m,-j},\underline{\pi_t}} \big\{ R^{m,j}_t(X^{l,m}_t,A^{l,m}_t) +  \delta V^{m,j}_{t+1} (\underline{F}(\underline{\pi_t},\tgamma^{l,m}_t,A_t^{l,m}), X_{t+1}^{m,j}) \big\lvert \underline{\pi_t},\hat{a}_{1:t-1}^{l,m},\hat{x}_{1:t}^{m,j} \big\} \nn\\
&> V^{m,j}_t(\underline{\pi_t}, \hat{x}_{t}^{m,j}).\label{b_eq:E8}
}
We will show that this leads to a contradiction.
Construct 
\begin{equation}
\hat{\gamma}^{m,j}_t(a_t^{m,j}|x_t^{m,j}) = \lb{\hat{\sigma}_t^{m,j}(a_t^{m,j}|\hat{a}_{1:t-1}^{l,m},\hat{x}_{1:t}^{m,j}) \;\;\;\;\; x_t^{m,j} = \hat{x}_t^{m,j} \\ \text{arbitrary} \;\;\;\;\;\;\;\;\;\;\;\;\;\; \text{otherwise.}  }
\end{equation}

Then for $\hat{a}_{1:t-1},\hat{x}_{1:t}^{m,j}$ and $\underline{\hat{\pi}_t}(x_t^{l,m}) = P^{\sigma}(x_t^{l,m}|\hat{z}_{1:t},\hat{a}_{1:t-1})$, we have
\seq{
\eq{
&V^{m,j}_t(\underline{\pi}_t, \hat{x}_t^{m,j})
 = \max_{\gamma_t^{m,j}(\cdot|\hat{x}_t^{m,j})} \E^{\tsigma^{l},\tsigma^{m},\gamma^{m,j}_t(\cdot|\hat{x}_t^{m,j}),\underline{\hat{\pi}_t} }\big\{ R^{m,j}_t(x^{l}_t,x_t^{m,-j},\hat{x}_t^{m,j},A_t^{l,m})\nn \\
&+ \delta V^{m,j}_{t+1} (\underline{F}(\underline{\hat{\pi}_t}  ,\tgamma^{l,m}_t,A^{l,m}_t), X_{t+1}^{m,j}) \big\lvert \underline{\hat{\pi}_t}, \hat{x}_{t}^{m,j} \big\}, \label{b_eq:E11}\\
\nn\\
&\geq\E^{\tsigma^{l},\tsigma^{m},\hat{\gamma}_t^{m,j}(\cdot|\hat{x}_t^{m,j}),\underline{\hat{\pi}_t}} \big\{ R^{m,j}_t(x^{l}_t,x_t^{m,-j},\hat{x}_t^{m,j},A_t^{l,m}) + \nn\\
&\delta V^{m,j}_{t+1} (\underline{F}(\underline{\hat{\pi}_t},\tgamma^{l,m}_t,A_t^{l,m}), {X}_{t+1}^{m,j}) \big\lvert \underline{\hat{\pi}_t},\hat{x}_{t}^{m,j} \big\}   
\\ \nn\\
&=\sum_{x^{l}_t,x^{m,-j}_t,a_t^{l,m},x_{t+1}^{m,j}}   \big\{ R^{m,j}_t(x^{l}_t,x_t^{m,-j},\hat{x}_t^{m,j},a_t^{l,m}) + \delta V^{m,j}_{t+1} (\underline{F}(\underline{\hat{\pi}_t},\tgamma^{l,m}_t,a^{l,m}_t), x_{t+1}^{m,j})\big\}
\nn\\
&\underline{\hat{\pi}_t}(x^{l,m}_t)\gamma_t^{l}(a_t^l|x_t^l)\gamma_t^{m,-j}(a_t^{m,-j}|x_t^{m,-j})\hat{\gamma}^{m,j}_t(a^{m,j}_t|\hat{x}_t^{m,j})
\\ \nn\\
&= \sum_{x^{l}_t,x_t^{m,-j},a_t^{l,m},x_{t+1}^{m,j}}  \big\{ R^{m,j}_t(x^{l}_t,x_t^{m,-j},\hat{x}_t^{m,j},a_t^{l,m}) + \delta V^{m,j}_{t+1} (\underline{F}(\underline{\hat{\pi}_t},\tgamma^{l,m}_t,a^{l,m}_t), x_{t+1}^{m,j})\big\}
\nn\\
&\underline{\hat{\pi}_t}(x^{l,m}_t)\gamma_t^{l}(a_t^l|x_t^l)\gamma_t^{m,-j}(a_t^{m,-j}|x_t^{m,-j})\hat{\sigma}^{m,j}_t(a_t^{m,j}|\hat{a}_{1:t-1}^{l,m},\hat{x}_{1:t}^{m,j})\label{b_eq:E9}\\
\nn\\
&= \E^{\tsigma^{l},\tsigma^{m},\underline{\hat{\pi}_t} } \big\{ R^{m,j}_t(x^{l}_t,x_t^{m,-j},\hat{x}_t^{m,j},A_t^{l,m})+\nn\\
&\delta V^{m,j}_{t+1} (\underline{F}(\underline{\hat{\pi}_t},\tgamma^{l,m}_t,A_t^{l,m}), X_{t+1}^{m,j}) \big\lvert \underline{\hat{\pi}_t},\hat{a}_{1:t-1}^{l,m}, \hat{x}_{1:t}^{m,j} \big\}  \\
&> V^{m,j}_t(\underline{\hat{\pi}_t}, \hat{x}_{t}^{m,j}), \label{b_eq:E10}
}
where \eqref{b_eq:E11} follows from definition of $V^{m,j}_t$ in \eqref{eq:Vdef}, \eqref{b_eq:E9} follows from definition of $\hat{\gamma}_t^{m,j}$ and \eqref{b_eq:E10} follows from \eqref{b_eq:E8}. However this leads to a contradiction.
}

\end{proof}

\begin{lemma}
\label{b_lemma:1}
Let $\tsigma$ be the strategies computed by the methodology in Section~III. 
Then $\forall t\in [T], a_{1:t-1}^{l,m},x_{1:t}^{m,j}$,
\begin{gather}
V^{m,j}_t(\underline{\pi_t}, x_t^{m,j}) =
\E^{\tsigma_{t:T}^{l,m},\underline{\pi_t}} \big\{ \sum_{n=t}^T \delta^{n-t}R^{m,j}_n(X_t^{l,m},A_t^{l,m}) \big\lvert \underline{\pi_t},  a_{1:t-1}^{l,m},x_{1:t}^{m,j} \big\} .
\end{gather} 
\end{lemma}

\begin{proof}
%
\seq{
We prove the lemma by induction. For $t=T$,
\eq{
 &\E^{\tsigma_t^{l},\tsigma_t^{m},\underline{\pi_t} } \big\{  R^{m,j}(X_T^{l,m},A_T^{l,m}) \big\lvert\underline{\pi_t}, a_{1:T-1}^{l,m},x_{1:T}^{m,j} \big\}\nn\\
 &= \sum_{x_{T}^{l,\{m,-j\}},a_T^{l,m}} R^{m,j}_T(x^{l,m}_T,a^{l,m}_T) \underline{\pi_t}(x^{l}_T,x^{m,-j}_T)\tsigma_{T}^{l,m}(a_T^{l,m}|\underline{\pi}_{T},x_{T}^{m,j}) \\
 &= V^{m,j}_T(\underline{\pi_t}, x_T^{m,j}) \label{b_eq:C1},
}
}
where \eqref{b_eq:C1} follows from the definition of $V^{m,j}_t$ in \eqref{eq:Vdef}.
Suppose the claim is true for $t+1$, i.e., $\forall  t\in [T],  a_{1:t}^{l,m}$,  $x_{1:t+1}^{m,j}$
\begin{gather}
V^{m,j}_{t+1}(\pi_{t+1},x_{t+1}^{m,j}) = \E^{\tsigma_{t+1:T}^{l,m},\pi_{t+1}}
\big\{ \sum_{n=t+1}^T \delta^{n-t-1}R^{m,j}_n(X^{l,m}_n,A^{l,m}_n) \big\lvert \pi_{t+1}, a_{1:t}^{l,m}, x_{1:t+1}^{m,j} \big\} 
\label{b_eq:CIndHyp}.
\end{gather}
Then $\forall  t\in [T],a_{1:t-1}^{l,m}, x_{1:t}^{m,j}$, we have
\seq{
\eq{
&\E^{\tsigma_{t:T}^{l,m},\underline{\pi_t} } \big\{ \sum_{n=t}^T \delta^{n-t} R^{m,j}_n(X^{l,m}_n,A^{l,m}_n) \big\lvert \underline{\pi_t},  a_{1:t-1}^{l,m},x_{1:t}^{m,j} \big\} 
\nonumber 
\\
&=  \E^{\tsigma_{t:T}^{l,m},\underline{\pi_t}} \big\{R^{m,j}_t(X^{l,m}_t,A^{l,m}_t)+\delta \E^{\tsigma_{t:T}^{l,m},\underline{\pi_t} } 
\nonumber \\ 
& \big\{ \sum_{n=t+1}^T \delta^{n-t-1}R^{m,j}_n(X^{l,m}_n,A^{l,m}_n)\big\lvert \underline{F}(\underline{\pi_t},\tgamma_t^{l,m},A^{l,m}_t), x_{1:t}^{m,j},X_{t+1}^{m,j}\big\} \big\lvert \underline{\pi_t},  a_{1:t-1}^{l,m}, x_{1:t}^{m,j} \big\} \label{b_eq:C2}
\\
&=  \E^{\tsigma_{t:T}^{l,m},\underline{\pi_t}} \big\{R^{m,j}_t(X_t,A_t) +\delta\E^{\tsigma_{t+1:T}^{l,m},\tsigma_{t+1:T}^{m,j},\underline{F}(\underline{\pi_t},\tgamma^{l,m}_t,A_t^{l,m})}
\nonumber 
\\
&\big\{ \sum_{n=t+1}^T \delta^{n-t-1}R^{m,j}_n(X^{l,m}_n,A^{l,m}_n)\big\lvert \underline{F}(\underline{\pi_t},\tgamma^{l,m}_t,A^{l,m}_t), x_{1:t}^{m,j},X_{t+1}^{m,j}\big\} \big\lvert \underline{\pi_t},  a_{1:t-1}^{l,m},x_{1:t}^{m,j} \big\} \label{b_eq:C3}
\\
&=  \E^{\tsigma_{t:T}^{l,m},\underline{\pi_t}} \big\{R^{m,j}_t(X^{l,m}_t,A^{l,m}_t) + \delta V^{m,j}_{t+1}(\underline{F}(\underline{\pi_t},\tgamma^{l,m}_t,A^{l,m}_t), X_{t+1}^{m,j}) \big\lvert \underline{\pi_t},  a_{1:t-1}^{l,m}, x_{1:t}^{m,j} \big\} 
\label{b_eq:C4}
\\
&=  \E^{\tsigma_t^{l},\tsigma_t^{m},\underline{\pi_t}} \big\{R^{m,j}_t(X^{l,m}_t,A^{l,m}_t) +  \delta V^{m,j}_{t+1}(\underline{F}(\underline{\pi_t},\tgamma^{l,m}_t,A^{l,m}_t), X_{t+1}^{m,j}) \big\lvert \underline{\pi_t},  a_{1:t-1}^{l,m}, x_{1:t}^{m,j} \big\} 
\label{b_eq:C5}
\\
&=V^{m,j}_{t}(\underline{\pi_t}, x_t^{m,j}) \label{b_eq:C6},
}
} where \eqref{b_eq:C3} follows from Lemma~\ref{b_lemma:3}, 
\eqref{b_eq:C4} follows from the induction hypothesis in \eqref{b_eq:CIndHyp} and \eqref{b_eq:C6} follows from the definition of $V^{m,j}_t$ in \eqref{eq:Vdef}.
\end{proof}

\begin{lemma}
\label{b_lemma:3}
$\forall  t\in \mathcal{T}, (a_{1:t}^{l,m}, x_{1:t+1}^{m,j})$ and
$\sigma^{m,j}_{t}$ 
\eq{
&\E^{\tsigma_{t:T}^{l},\sigma_{t:T}^{m,j},\tsigma_{t:T}^{m,-j},\,\underline{\pi_t}}  \big\{ \sum_{n=t+1}^T R_n^{m,j}(X^{l,m}_n,A^{l,m}_n) \big\lvert \pi_{t}, a_{1:t}^{l,m}, x_{1:t+1}^{m,j} \big\} =\nn\\
& \E^{\tsigma^{l}_{t+1:T},\sigma_{t+1:T}^{m,j},\tsigma_{t+1:T}^{m,-j},\underline{F}(\underline{\pi_t},\tgamma^{l,m}_t,A^{l,m}_t)}  \big\{ \sum_{n=t+1}^T R_n^{m,j}(X^{l,m}_n,A^{l,m}_n) \big\lvert \pi_{t+1}, a_{1:t}^{l,m}, x_{1:t+1}^{m,j} \big\}. \label{eq:1}
}

\end{lemma}
\begin{proof} 
Since the above expectations involve random variables $X_{t+1}^{l,\{m,-j\}},A^{l,m}_{t+1:T}$, \\$X^{l,m}_{t+2:T}$, we consider the probability 
\seq{
\eq{
&\p^{\tsigma^{l}_{t:T},\sigma^{m,j}_{t:T},\tsigma_{t:T}^{m,-j}\,\underline{\pi_t}} (x_{t+1}^{l,\{m,-j\}},a^{l,m}_{t+1:T}, x^{l,m}_{t+2:T}\big\lvert \underline{\pi_t},  a_{1:t}^{l,m}, x_{1:t+1}^{m,j} ) = \frac{Nr}{Dr} \label{b_eq:2}
}
\vspace{-0.8cm}
\eq{
&\text{where    } Nr= \nn\\
&\sum_{\substack{x^{l, \{m,-j\}}_t,\\a_t^{l,m}}}\p^{\tsigma^{l}_{t:T},\sigma^{m,j}_{t:T},\tsigma_{t:T}^{m,-j}\,\underline{\pi_t}} (x_t^{l,\{m,-j\}},a^{l,m}_t, x^{l,\{m,-j\}}_{t+1},a^{l,m}_{t+1:T},x^{l,m}_{t+2:T} \big\lvert\underline{\pi_t}, a_{1:t-1}^{l,m}, x_{1:t}^{m,j} ) \\
&= \sum_{\substack{x_t^{l,\{m,-j\}},\\a_t^{l,m}}}\p^{\tsigma^{l}_{t:T},\sigma^{m,j}_{t:T},\tsigma_{t:T}^{m,-j}\,\underline{\pi_t}} (x^{l,\{m,-j\}}_t \big\lvert \underline{\pi_t}, a^{l,m}_{t+1:T}, x_{1:t}^{m,j} )\tsigma_t^{l,m}(a_t^{l,m}|\underline{\pi_t},x_t^{l,m})
\nonumber 
\\
&Q^{l,\{m,-j\}}(x^{l,\{m,-j\}}_{t+1}|x^{l,m}_t, a^{l,m}_t)\p^{\tsigma^{l,m}_{t:T}, \,\underline{\pi_t}} ( a^{l,m}_{t+1:T},x^{l,m}_{t+2:T}| a_{1:t-1}^{l,m} ,x_{1:t-1}^{m,j}, x^{l,m}_{t:t+1}) 
\\
=&\sum_{x_t^{l,\{m,-j\}}}\underline{\pi_t}(x_t^{l,\{m,-j\}})) \tsigma_t^{l,m}(a_t^{l,m}|\underline{\pi_t},x_t^{l,m}) Q^{l,\{m,-j\}}(x^{l,\{m,-j\}}_{t+1}|x^{l,m}_t, a^{l,m}_t)
\\
&\p^{\tsigma^{l,m}_{t+1:T} ,\, \pi_{t+1}} ( a^{l,m}_{t+1:T},x^{l,m}_{t+2:T}| \underline{\pi_t},a_{1:t-1}^{l,m},x_{1:t-1}^{m,j},x^{l,m}_{t:t+1}),\label{b_eq:Nr2}
}
where \eqref{b_eq:Nr2} follows from the conditional independence of types given common information, as shown in Claim~1 in~Appendix~\ref{app:0}, and the fact that probability on $(a^{l,m}_{t+1:T},x^{l,m}_{2+t:T})$ given $ x_{1:t}^{l,m},x^{l,m}_{t+1}, \pi_{t} $ depends on $ a_{1:t-1}^{l,m},x_{1:t}^{l,m},x^{l,m}_{t+1}, \pi_{t+1} $ through ${ \tsigma_{t+1:T}^{l,m} }$. Similarly, the denominator in \eqref{b_eq:2} is given by
\eq{
Dr &= \sum_{\tilde{x}_{t}^{l,\{m,-j\}},a_t^{l,m}} \p^{ \tsigma^{l}_{t:T},\sigma^{m,j}_{t:T},\tsigma_{t:T}^{m,-j},\,\underline{\pi_t}} (\tilde{x}_t^{l,\{m,-j\}}, a^{l,\{m,-j\}}_t,x_{t+1}^{m,j}\big\lvert \underline{\pi_t},a_{1:t-1}^{l,m}, x_{1:t}^{m,j} )\\
=&\sum_{\tilde{x}_{t}^{l,m\{m,-j\}},a_t^{l,\{m,-j\}}} \underline{\pi_t}(\tilde{x}_t^{l,\{m,-j\}} \sigma_t^{m,j}(a_t^{m,j}| a_{1:t-1}^{l,m},x_{1:t}^{l,m})  \nonumber\\
&\tsigma_t^{l,\{m,-j\}}(a_t^{l,\{m,-j\}}|\underline{\pi_t}, \tilde{x}_t^{l,\{m,-j\}}) Q^{l,\{m,j\}}(x^{\{m,j\}}_{t+1}|x^{l,m}_t, a^{l,m}_t)\label{b_eq:4}
%
}

By canceling the terms $\sigma_t^{m,j}(\cdot)$ and $Q^{m,j}(\cdot)$ in the numerator and the denominator, \eqref{b_eq:2} is given by
\eq{
&\hspace{-23pt}\frac{\sum_{x_t^{l,\{m,-j\}}}\underline{\pi_t}(x_t^{l,\{m,-j\}}) \tsigma_t^{l,\{m,-j\}}(a_t^{l,\{m,-j\}}|\underline{\pi_t}, x_t^{l,\{m,-j\}}) Q_{t+1}^{l,\{m,-j\}}(x^{l,m\{m,-j\}}_{t+1}|x^{l,m}_t, a^{l,m}_t)}{\sum_{\tilde{x}_{t}^{l,\{m,-j\}}} \underline{\pi_t}(\tilde{x}_t^{l,\{m,-j\}}) \tsigma_t^{l,\{m,-j\}}(a_t^{l,\{m,-j\}}|\underline{\pi_t},\tilde{x}_t^{l,m})}  \nonumber \\
&\times\p^{\tsigma^{l,m}_{t+1:T}, \sigma^{m,j}_{t+1:T} ,\, \pi_{t+1}} ( a^{l,m}_{t+1:T},x^{l.m}_{t+2:T}|\underline{\pi_t}, a^{l,m}_{1:t-1},x_{1:t}^{m,j}, x^{l,m}_{t+1})\\
=&\underline{\pi}_{t+1}^{l,\{m,-j\}}(x_{t+1}^{l,\{m,-j\}}) \p^{ \sigma_{t+1:T}^{m,j},\tsigma^{l,\{m,-j\}}_{t+1:T},\tsigma^{m,j}_{t+1:T}\, \underline{\pi}_{t+1}} (a^{l,m}_{t+1:T},x^{l,m}_{t+2:T}|\underline{\pi_t} ,a_{1:t-1}^{l,m},x_{1:t}^{m,j}, x^{l,m}_{t+1})\label{b_eq:6}\\
=& \p^{\tsigma^{l}_{t+1:T},\sigma^{m,j}_{t+1:T},\tsigma_{t+1:T}^{m,-j}\,\underline{\pi}_{t+1} } (x_{t+1}^{l,\{m,-j\}},a_{1:t}^{l,m},x^{l,m}_{t+2:T} |\pi_{t}, a_{1:t-1}^{l,m}, x_{1:t+1}^{m,j} ),
}
}
where \eqref{b_eq:6} follows from using the definition of $\pi^{l,\{m,-j\}}_{t+1}(x_{t+1}^{l,\{m,-j\}})$ in \eqref{eq:piupdate}.

\end{proof}

\section{Part 2: Leader}
\label{app:P2}
In the following, we will show that, $\forall t, a_{1:t-1}^{l,m}, x_{1:t}^{l}, \sigma^{l}$
\eq{
&\E^{\tsigma^{l},\tsigma^{m},\underline{\pi_t}} \big\{ \sum_{n=t}^T \delta^{n-t}R_n^{l}(X_n^{l,m},A_n^{l,m}) |\underline{\pi_t},a_{1:t-1}^{l,m},x_{1:t}^{l}\big\} \nn\\
&\geq
 \E^{\sigma^{l},\hat{\sigma}^{m},\underline{\pi_t}} \big\{ \sum_{n=t}^T \delta^{n-t}R_n^{l}(X_n^{l,m},A_n^{l,m}) |\underline{\pi_t},a_{1:t-1}^{l,m},x_{1:t}^{l}\big\},
}
where $\hat{\sigma}^{m}\in BR^{m}(\sigma^{l},\hat{\sigma}^{m})$.

\begin{proof}
We prove the above result using induction and from results in Lemma~\ref{l_lemma:2} and \ref{l_lemma:1} proved in Appendix~\ref{l_app:lemmas}. 

For the base case at $t=T$, $\forall a_{1:T-1}^{l,m}, x_{1:T}^{l}, \sigma^{l}$
\seq{
\eq{
&\E^{\tsigma_{T}^{l},\tsigma_T^m, \underline{\pi_t}}\big\{  R_T^{l}(X^{l,m}_T,A^{l,m}_T) \big\lvert \underline{\pi_T}, a_{1:T-1}^{l,m},x_{1:T}^{l}\big\}\nn\\
&=V^{l}_T(\underline{\pi_T},x^{l}_T)  \label{l_eq:T2a}\\
&\geq \E^{ \sigma_{T}^{l}(\cdot|x_T),\hat{\sigma}_T^m,\underline{\pi}_T}\big\{ R_T^{l}(X^{l,m}_T,A^{l,m}_T) \big\lvert \underline{\pi_T}, a_{1:T-1}^{l,m},x_{1:T}^{l} \big\} \label{l_eq:T2},\nn\\
&\text{ where } \hat{\sigma}_T^{m}\in BR_T^{m}(\underline{\pi_T},a_{1:T-1}^{l,m},\sigma_T^{l},\hat{\sigma}_T^{m})
}
}
where (\ref{l_eq:T2a}) follows from Lemma~\ref{l_lemma:1} and (\ref{l_eq:T2}) follows from Lemma~\ref{l_lemma:2} in Appendix~\ref{l_app:lemmas}. Let the induction hypothesis be that for $t+1$, $\forall  a_{1:t}^{l,m},x_{1:t+1}^{l,m}, \sigma^{l}$,
\seq{
\eq{
 & \E^{ \tsigma_{t+1:T}^{l},\tsigma^m_{t:T},\pi_{t+1}} \big\{ \sum_{n=t+1}^T R_n^{l}(X_n^{l,m},A_n^{l,m}) \big\lvert \pi_{t+1}, a_{1:t}^{l,m},x_{1:t+1}^{l} \big\} \nn\\
  &\geq \E^{ \sigma_{t+1:T}^{l},\hat{\sigma}_{t:T}^m,\pi_{t+1}} \big\{ \sum_{n=t+1}^T R_n^{l}(X_n^{l,m},A_n^{l,m}) \big\lvert \pi_{t+1}, a_{1:t}^{l,m},x_{1:t+1}^{l} \big\} \label{l_eq:PropIndHyp}\\
  &\text{where } 
  \hat{\sigma}^{m}_{t+1:T}\in BR_{t+1}^{m}(\pi_{t+1},a_{1:t}^{l,m},\sigma_{t+1:T}^{l},\hat{\sigma}_{t+1:T}^{m})
}
Then $\forall a_{1:t-1}^{l,m},x_{1:t}^{l}, \sigma^{l}$, we have
\eq{
&\E^{ \tsigma_{t:T}^{l,m},\underline{\pi_t}} \big\{ \sum_{n=t}^T R_n^{l}(X_n^{l,m},A_n^{l,m}) \big\lvert \pi_{t}, a_{1:t-1}^{l,m},x_{1:t}^{l}\big\} \nn \\
&= V^{l}_t(\underline{\pi_t},x^{l}_t)\label{l_eq:T1}\\
&\geq \E^{\gamma_t^{l},\hat{\gamma}_t^m,\underline{\pi_t}} \big\{ R_t^{l}(X^{l,m}_t,A^{l,m}_t) + V_{t+1}^{l} (\underline{F}(\underline{\pi_t},\tgamma_t^{l},\hat{\gamma}_t^m,A_t^{l,m}),X_{t+1}^{l})\big\vert \pi_{t}, a_{1:t-1}^{l,m},x_{1:t}^{l} \big\}  \label{l_eq:T3}\\
&= \E^{ \sigma_t^{l},\hat{\sigma}_t^m,\underline{\pi_t}} \big\{ R_t^{l}(X^{l,m}_t,A^{l,m}_t) +\nn\\
&\E^{\tsigma_{t+1:T}^{l,m},\underline{F}(\underline{\pi_t},\tgamma_t^{l},\hat{\gamma}_t^m,A_t^{l,m})}  \big\{ \sum_{n=t+1}^T R_n^{l}(X^{l,m}_n,A^{l,m}_n)  \big\lvert x_{1:t}^{l},X_{t+1}^{l}\big\}  \big\vert \pi_{t}, a_{1:t-1}^{l,m},x_{1:t}^{l}\big\}  \label{l_eq:T3b}\\
&\geq \E^{\sigma_t^{l},\hat{\sigma}_t^m,\underline{\pi_t}} \big\{ R_t^{l}(X^{l,m}_t,A^{l,m}_t)+\E^{\sigma_{t+1:T}^{l},\hat{\sigma}_{t+1:T}^{m}, \underline{F}(\underline{\pi_t},\tgamma_t^{l},\hat{\gamma}_t^m,A_t^{l,m})}  
}
\eq{
& \big\{ \sum_{n=t+1}^T R_n^{l}(X_n^{l,m},A_n^{l,m})  \big\lvert F(\underline{\pi_t},\tgamma_t^{l},\hat{\gamma}_t^m),x_{1:t}^{l},X_{t+1}^{l}\big\} \big\vert \pi_{t}, a_{1:t-1}^{l,m},x_{1:t}^{l} \big\}  \label{l_eq:T4} \\
&= \E^{{\sigma}_t^{l},\hat{\sigma}_t^m,\underline{\pi_t}} \big\{ R_t^{l}(X^{l,m}_t,A^{l,m}_t) +\nn\\
&\E^{\sigma_{t:T}^{l},\hat{\sigma}_{t:T}^m,\underline{\pi_t}}\big\{ \sum_{n=t+1}^T R_n^{l}(X_n^{l,m},A_n^{l,m})  \big\lvert (\underline{\pi_t},\gamma_t^{l},\hat{\gamma}_t^m),x_{1:t}^{l},X_{t+1}^{l}\big\} \big\vert \pi_{t}, a_{1:t-1}^{l,m},x_{1:t}^{l}\big\}  \label{l_eq:T5}\\
&=\E^{\sigma_{t:T}^{l},\hat{\sigma}_{t:T}^m,\underline{\pi_t}} \big\{ \sum_{n=t}^T R_n^{l}(X_n^{l,m},A_n^{l,m}) \big\lvert \pi_{t}, a_{1:t-1}^{l,m},x_{1:t}^{l}\big\}  \label{l_eq:T6},
}
}
where  $\hat{\sigma}_t^m \in BR_t^m(\underline{\pi_t},a_{1:t-1}^{l,m},\sigma_t^{l},\hat{\sigma}_t^m)$,  (\ref{l_eq:T1}) follows from Lemma~\ref{l_lemma:1}, (\ref{l_eq:T3}) follows from Lemma~\ref{l_lemma:2}, (\ref{l_eq:T3b}) follows from Lemma~\ref{l_lemma:1} and 
(\ref{l_eq:T4}) follows from induction hypothesis in (\ref{l_eq:PropIndHyp}), (\ref{l_eq:T5}) follows from the fact that probability on $(a^{l,m}_{t+1:T},x^{l,m}_{2+t:T})$ given $\underline{\pi_t}, a_{1:t}^{l,m},x_{1:t+1}^{l,m}$ depends on $\pi_{t+1}, a_{1:t}^{l,m},x_{1:t+1}^{l,m}$ through ${ \tsigma_{t+1:T}^{l,m} }$. 
\end{proof}

\section{}
\label{l_app:lemmas}
\begin{lemma}
\label{l_lemma:2}
$\forall t\in [T], \pi_t,a_{1:t-1}^{l,m},x_{1:t}^{l}, \forall \sigma^{l}$
\eq{
&V_t^{l}(\underline{\pi_t},x^{l}_t) \geq \nn\\
&\E^{\sigma_t^{l},\bar{\sigma}_t^m,\underline{\pi_t}} \big\{ R_t^{l}(X^{l,m}_t,A^{l,m}_t) + V_{t+1}^{l} (\underline{F}(\underline{\pi}_{t},\tgamma_t^l,\bar{\gamma}_t^m,A_t^{l,m}),X_{t+1}^{l}) \big\lvert \underline{\pi}_{t}, a_{1:t-1}^{l,m},x_{1:t}^{l}\big\}\label{l_eq:lemma2}
}
where 
\eq{
%
%
\forall j, \bar{\sigma}_t^{m,j}&\in\arg\max_{\sigma_t^{m,j}} \bigcap_{a_{1:t-1}^{l,m},x_{1:t}^{m,j}}BR^{m,j}_t(\underline{\pi_t},a_{1:t-1}^{l,m},x_{1:t}^{m,j},\sigma_t^{l},\sigma_t^{m,j},\bar{\sigma}_{t}^{m,-j},\tsigma_{t+1:T}^{l,m})\\ 
{\tgamma}^{m}_t &= {\tsigma}_t^{m}(\cdot|a_{1:t-1}^{l,m},\cdot)\\
\gamma_t^l &= \sigma_t^l(\cdot|a^{l,m}_{1:t-1},x^{l}_{1:t-1},\cdot)\\
\bar{\gamma}_t^m &\in \bar{BR}_t^m(\underline{\pi_t},\gamma_t^{l},\bar{\gamma}_t^m)
}
where we assume that $\bar{\sigma}^m$ are of type $m$ (Since if they are not, one can find an equivalent policies of type $m$ that achieve same reward profile, as shown in Appendix~\ref{app:gsm}).

\end{lemma}

\begin{proof}

We prove this by contradiction. Suppose the claim is not true for $t$. This implies $\exists\  \breve{\sigma}_t^{l},\hat{a}_{1:t-1}^{l,m}$, $\hat{x}_{1:t}^{l}$ such that $\hat{\pi}_t = P^{\tsigma^{l,m}_{1:t}}(\cdot|\hat{a}_{1:t-1}^{l,m},\hat{x}_{1:t-1}^{l},\cdot)$ and,
\eq{
&\E^{ \breve{\sigma}_t^{l},\hat{\sigma}_t^m,\underline{\hat{\pi}_t}} \big\{ R_t^{l}(X^{l,m}_t,A^{l,m}_t) + V_{t+1}^{l} (\underline{F}(\underline{\hat{\pi}_t},{\tgamma}_t^{l},\hat{\gamma}_t^m, A^{l,m}_t),X_{t+1}^{l}) \big\lvert \hat{\underline{\pi}}_{t}, \hat{z}_{1:t},\hat{a}_{1:t-1}^{l,m},\hat{x}_{1:t}^{l}\big\} >\nn\\ &V_t^{l}({\hat{\underline{\pi}}_t},\hat{x}_t^{l}),\label{l_eq:E8}
}
where $\breve{\gamma}^{l}_t(\cdot) = \breve{\sigma}^{l}_t(\cdot|\hat{a}_{1:t-1}^{l,m},x^{l}_{1:t-1},\cdot)$, $\hat{\sigma}_t^m\in BR_t^m(\underline{\hat{\pi}_t},\hat{a}_{1:t-1}^{l,m},\breve{\sigma}_t^{l},\hat{\sigma}_t^m,{\tsigma}_{t+1:T}^{l,m})$ and $\hat{\gamma}_t^{m}$ satisfies
\eq{
\hat{\gamma}^{m}_t &= \hat{\sigma}_t^{m}(\cdot|\hat{a}_{1:t-1}^{l,m},\hat{x}_{1:t-1}^{m},\cdot)
}


Then for $\hat{a}_{1:t-1}^{l,m},\hat{x}_{1:t}^{l}$, we have
\seq{
\eq{
&V_t^{l}(\underline{\pi_t},\hat{x}_t^{l})= \nn\\
&\max_{\gamma^{l}_t(\cdot|\hat{x}_t^l)} \E^{\gamma^{l}_t,\bar{\gamma}^{m}_t,\underline{\pi_t} } \big\{ R_t^{l}(X^{l,m}_t,A^{l,m}_t) + V_{t+1}^{l} (\underline{F}(\underline{\pi_t},\tgamma^{l}_t,\bar{\gamma}^{m}_t,A^{l,m}_t),X_{t+1}^{l}) \big\lvert \hat{\pi}_{t}, \hat{a}_{1:t-1}^{l,m},\hat{x}_{1:t}^{l}\big\} \label{l_eq:E11}\\
&\geq\E^{\breve{\gamma}^{l}_t,\hat{\gamma}^{m}_t,\underline{\pi_t}} \big\{ R_t^{l}(X^{l,m}_t,A^{l,m}_t)+V_{t+1}^{l} (\underline{F}(\underline{\pi_t},{\tgamma}^{l}_t,\hat{\gamma}^{m}_t,A^{l,m}_t),X_{t+1}^{l}) \big\lvert \hat{\pi}_{t}, \hat{a}_{1:t-1}^{l,m},\hat{x}_{1:t}^{l} \big\}   \\
&= \E^{\breve{\sigma}^{l}_t,\hat{\sigma}^{m}_t,\underline{\pi_t}} \big\{ R_t^{l}(X^{l,m}_t,A^{l,m}_t)+  V_{t+1}^{l} (\underline{F}(\underline{\hat{\pi}_t},{\tgamma}^{l}_t,\hat{\gamma}^{m}_t,A^{l,m}_t),X_{t+1}^{l}) \big\lvert \hat{\underline{\pi}}_{t}, \hat{a}_{1:t-1}^{l,m},\hat{x}_{1:t}^{l}\big\}  \label{l_eq:E9}\\
&> V_t^{l}(\underline{\hat{\pi}_t},\hat{x}_t^{l}) \label{l_eq:E10} 
}
where (\ref{l_eq:E11}) follows from the definition of $V_t^{l}$ in (\ref{eq:Vdef}), (\ref{l_eq:E9}) follows from definition of $\breve{\gamma}_t^{l},\hat{\gamma}^{m}_t$ and (\ref{l_eq:E10}) follows from (\ref{l_eq:E8}). However this leads to a contradiction. 
}
\end{proof}

\begin{lemma}
\label{l_lemma:1}
$\forall t\in [T],a_{1:t-1}^{l,m},x_{1:t}^{l}$
\eq{
V^{l}_t(\underline{\pi_t},x^{l}_t)&= \E^{\tsigma_{t:T}^{l,m} ,\underline{\pi_t}} \big\{ \sum_{n=t}^T R_n^{l}(X_n^{l,m},A_n^{l,m}) \big\lvert \underline{\pi}_{t}, a_{1:t-1}^{l,m},x_{1:t}^{l}\big\} .
}
\end{lemma}
\begin{proof}
\seq{
We prove the lemma by induction. For $t=T$, 
\eq{
 &\E^{\tsigma_{T}^{l,m},\underline{\pi_T}} \big\{  R_T^{l}(X^{l,m}_T,A^{l,m}_T) \big\lvert \pi_{T}, a_{1:T-1}^{l,m},x_{1:T}^{l}\big\}\nn \\
 &= \sum_{x_T^{m}, a^{l,m}_T} \underline{\pi}_T(x_T^{m})R_T^{l}(x^{l,m}_T,a^{l,m}_T) \tsigma_{T}^{l,m}(a^{l,m}_T|\underline{\pi}_T,x_T^{l,m})\\ 
 &=V^{l}_T(\underline{\pi_T},x^{l}_T) \label{l_eq:C1}
}
}
where (\ref{l_eq:C1}) follows from the definition of $V_t^{l}$ in (\ref{eq:Vdef}).

Suppose the claim is true for $t+1$, i.e., $\forall  t\in [T], a_{1:t}^{l,m},x_{1:t+1}^{l}$
\eq{
&V^{l}_{t+1}(\pi_{t+1},x_{t+1}^{l}) = \E^{ \tsigma_{t+1:T}^{l,m},\pi_{t+1}} \big\{ \sum_{n=t+1}^T R_n^{l}(X_n^{l,m},A_n^{l,m}) \big\lvert\pi_{t+1}, a_{1:t}^{l,m},x_{1:t+1}^{l}\big\} \label{l_eq:CIndHyp}.
}

Then $\forall  t\in [T], a_{1:t-1}^{l,m},x_{1:t}^{l}$, we have
	\seq{
\eq{
&\E^{ \tsigma_{t:T}^{l,m},\underline{\pi_t}} \big\{ \sum_{n=t}^T R_n^{l}(X_n^{l,m},A_n^{l,m}) \big\lvert \underline{\pi_t}, a_{1:t-1}^{l,m},x_{1:t}^{l} \big\} \nn \\
&=  \E^{ \tsigma_{t:T}^{l,m}, \underline{\pi_t} } \big\{R_t^{l}(X^{l,m}_T,A^{l,m}_t) +\nn\\
&\E^{ \tsigma_{t:T}^{l,m},\underline{\pi_t} } \big\{ \sum_{n=t+1}^T R_n^{l}(X_n^{l,m},A_n^{l,m}) \big\lvert a_{1:t-1}^{l,m},A^{l,m}_t,x_{1:t}^{l},X_{t+1}^{l}\big\} \big\lvert \underline{\pi_t}, a_{1:t-1}^{l,m},x_{1:t}^{l}\big\} \label{l_eq:C2}\\
&=  \E^{\tsigma_{t:T}^{l,m},\underline{\pi_t}} \big\{R_t^{l}(X^{l,m}_T,A^{l,m}_t)+\E^{\tsigma_{t+1:T}^{l,m} ,\underline{F}(\underline{\pi_t},\tgamma^{l,m}_t,A^{l,m}_t)}\big\{ \sum_{n=t+1}^T R_n^{l}(X_n^{l,m},A_n^{l,m}) \big\lvert\nn\\
&  a_{1:t-1}^{l,m},A^{l,m}_t,x_{1:t}^{l},X_{t+1}^{l}\big\} \big\lvert \underline{\pi_t}, a_{1:t-1}^{l,m},x_{1:t}^{l}\big\} \label{l_eq:C3}\\
&=  \E^{\tsigma_{t}^{l,m},\underline{\pi_t}} \big\{R_t^{l}(X^{l,m}_T,A^{l,m}_t) + V^{l}_{t+1}(\underline{F}(\underline{\pi_t},\tgamma^{l,m}_t,A^{l,m}_t),X_{t+1}^{l}) \big\lvert a_{1:t-1}^{l,m},x_{1:t}^{l}\big\} \label{l_eq:C5}\\
&=V^{l}_{t}(\underline{\pi_t},x^{l}_t) \label{l_eq:C6},
}
}
where (\ref{l_eq:C3}) follows from~Lemma~\ref{l_lemma:3}, (\ref{l_eq:C5}) follows from the induction hypothesis in (\ref{l_eq:CIndHyp}), 
and (\ref{l_eq:C6}) follows from the definition of $V_t^{l}$ in (\ref{eq:Vdef}).
\end{proof}

\begin{lemma}
\label{l_lemma:3}
$\forall t\in \mathcal{T}, \sigma_t^{l}, (a_{1:t}^{l,m},x_{1:t+1}^{l})$ 
\eq{
&\E^{\sigma_{t:T}^{l},\tsigma_{t:T}^{m},\underline{\pi_t}}  \big\{ \sum_{n=t+1}^T R_n^{l}(X^{l,m}_n,A^{l,m}_n) \big\lvert \underline{\pi_t}, a_{1:t}^{l,m}, x_{1:t+1}^{l} \big\} =\nn\\
& \E^{\sigma_{t+1:T}^{l},\tsigma_{t+1:T}^{m},\underline{F}(\underline{\pi_t},\tgamma_t^{l,m},a_t^{l,m})}  \big\{ \sum_{n=t+1}^T R_n^{l}(X^{l,m}_n,A^{l,m}_n) \big\lvert \underline{\pi_t}, a_{1:t}^{l,m}, x_{1:t+1}^{l} \big\}, \label{l_eq:1}
}
where ${\tgamma}_t^{l,m}=\tsigma^{l,m}_t(\cdot|\underline{\pi}_t,\cdot)$.
\end{lemma}
\begin{proof} 
Since the above expectations involve random variables $X^{l,m}_{t+1:T},A^{l,m}_t $, we consider the probability 
\seq{
\eq{
&\p^{\sigma_{t:T}^{l},\tsigma_{t:T}^{m},\underline{\pi_t}} (x_{t+1}^{m},x_{t+2:T}^{l,m},a_{t+1:T}^{l,m}\big\lvert \underline{\pi_t},  a_{1:t}^{l,m}, x_{1:t+1}^{l} ) = \frac{Nr}{Dr} \label{l_eq:2}
}
\vspace{-0.4cm}
\eq{
&\text{where}\nn\\ 
&Nr 
=\sum_{x_t^{m},a_t^{l,m}}\p^{\sigma_{t:T}^{l},\tsigma_{t:T}^{m},\underline{\pi_t}} (x_t^{m},a^{l,m}_t, x^{m}_{t+1},a^{l,m}_{t+1:T},x^{l,m}_{t+2:T} \big\lvert\underline{\pi_t}, a_{1:t-1}^{l,m}, x_{1:t}^{l} ) \\
=&\sum_{x_t^{m},a_t^{l,m}}\underline{\pi}_t(x_t^{m})\sigma_t^{l,m}(a_t^{l,m}|a_{1:t-1}^{l,m}, x_{1:t}^{l,m})   Q^{l,m}(x^{m}_{t+1}|x^{l,m}_t, a_t^{l,m})
\\
&\p^{\sigma_{t+1:T}^{l},\tsigma_{t+1:T}^{m},\underline{\pi}_{t}} (x^{l,m}_{t+2:T}, a^{l,m}_{t+1:T}| \underline{\pi_t},a_{1:t}^{l,m},x_{1:t}^{l},x_{t+1}^{l,m}),\label{l_eq:Nr2}
}
where \eqref{l_eq:Nr2} follows from the fact that probability on $(x^{l,m}_{t+2:T},a^{l,m}_{t+1:T})$ given\\ $\underline{\pi_t},a_{1:t}^{l,m},x_{1:t}^{l},x_{t+1}^{l,m}$ depends on $\underline{\pi_t} ,a_{1:t}^{l,m},x_{1:t}^{l},x_{t+1}^{l,m}$ through $\sigma_{t+1:T}^{l},\tsigma_{t+1:T}^{m}$. Similarly, the denominator in \eqref{l_eq:2} is given by
\eq{
Dr &=  \sum_{x_t^{m}}\p^{ \sigma_{t:T}^{l},\tsigma_{t:T}^{m},\, \underline{\pi_t}} ( x_t^{m},a^{l,m}_t,x_{t+1}^{l}\big\lvert \underline{\pi_t}, a_{1:t-1}^{l,m}, x_{1:t}^{l} )\\
&=\sum_{x_t^{m}}\underline{\pi}_t(x_t^{m})\tsigma_t^{l,m}(a_t^{l,m}|\underline{\pi_t},x_t^{l,m})Q^{l}(x^{l}_{t+1}|x^{l,m}_t, a^{l,m}_t)\label{l_eq:3}
}

By canceling the terms $Q^{l}(\cdot)$ in the numerator and the denominator, \eqref{l_eq:2} is given by
\eq{
&\frac{\sum_{x_t^{l}}\underline{\pi}_t(x_t^{m}) \tsigma_t^{m}(a_t^{m}|\underline{\pi_t},x_t^{m}) Q_{t+1}^{m}(x^{m}_{t+1}|x^{l,m}_t,a^{l,m}_t)}{\sum_{\tilde{x}_{t}^{m}} \underline{\pi}_t(\tilde{x}_t^{m}) \tsigma_t^{m}(a_t^{m}|\underline{\pi_t},\tilde{x}_t^{m})}  \nonumber \\
&\times\p^{\tsigma^{l,m}_{t+1:T},\, \pi_{t+1}} ( a_{t+1:T}^{l,m},x^{l,m}_{t+2:T}| a_{1:t-1}^{l,m},x_{1:t}^{l}, x_{t+1})\\
=&\pi_{t+1}^{m}(x_{t+1}^{m}) \p^{\sigma_{t+1:T}^{l},\tsigma_{t+1:T}^{m} ,\, \pi_{t+1}} ( a_{t+1:T}^{l,m},x^{l,m}_{t+2:T}| a_{1:t-1}^{l,m},x_{1:t}^{l}, x_{t+1})\label{l_eq:6}\\
=& \p^{\sigma_{t+1:T}^{l},\tsigma_{t+1:T}^{m}\, \pi_{t+1} } (x_{t+1}^{m},x_{t+2:T}^{l,m},a_{t+1:T}^{l,m} | a_{1:t-1}^{l,m}, x_{1:t+1}^{l} ),
}

where \eqref{l_eq:6} follows from using the definition of $\underline{\pi}_{t+1}(x_{t+1}^{l,m})$ in \eqref{eq:piupdate}.
}
\end{proof}

\section{Extra Lemmas}
\begin{lemma}
\label{lemma:BR}
Let $\sigma_t^l$ be any strategy of the leader. Let
\eq{
\hat{\sigma}^{m} &\in BR^{m}(\tsigma_{1:t-1}^l,\sigma_t^{l},\tsigma_{t+1:T}^l,\hat{\sigma}^m)\label{b_eq:EL2}}
where we assume that $\hat{\sigma}^m$ are of type $m$ (Since if they are not, one can find an equivalent policies of type $m$ that achieve same reward profile, as shown in Appendix~\ref{app:gsm}).
Let
\eq{
\hat{\gamma}_t^m &\in \bar{BR}_t^m(\pi_t,\gamma_t^{l},\hat{\gamma}_t^m)\label{b_eq:EL4}
}
where $\gamma_t^{l} = \sigma_t^{l}(\cdot|a_{1:t-1}^{l,m},x_{1:t-1}^{l},\cdot)$. Then $\forall a_{1:t-1},x_{1:t-1}^{m,j},x_{1:t-1}^{m,j},x_{1:t-1}^{l}$, and for every that satisfy~\eqref{b_eq:EL4} $\exists (\hat{\sigma}_t^m)$ that satisfy~\eqref{b_eq:EL2} such that
\eq{
\hat{\gamma}_t^{m,j} &= \hat{\sigma}_t^{m,j}(\cdot|a^{l,m}_{1:t-1},\cdot)
}where,

\eq{
&BR_t^{m,j}(\underline{\pi_t},a_{1:t-1}^{l,m},x_{1:t}^{m,j},\sigma_{t:T}^{l},\sigma_{t:T}^{m,-j}) \nn\\
&:=\arg\max_{ \sigma^{m,j}} \E^{\sigma_{t:T}^{l}\sigma_{t:T}^{m,j}\sigma_{t:T}^{m,-j},\sigma_{t:T}^{m,j},\underline{\pi_t}} \big\{ \sum_{n=t}^T \delta^{n-t}R_n^{m,j}(X_n^{l,m},A_n^{l,m}) |\underline{\pi_t},a_{1:t-1}^{l,m},x_{1:t}^{m,j}\big\}\\
&BR^{m,j}(\sigma^l,\sigma^{m,-j},\sigma^{m,j}) :=\bigcap_t \bigcap_{a_{1:t-1}^{m,j}}\bigcap_{x_{1:t}^{m,j}}  BR_t^{m,j}(\underline{\pi_t},a_{1:t-1}^{l,m},x_{1:t}^{m,j},\sigma_{t:T}^{l},\sigma_{t:T}^{m,-j}).
}

where
\eq{
\pi_t(\cdot)=P^{\tsigma_{1:t-1}^{l,m}}(\cdot|a_{1:t-1}). \label{eq:pi_def}
}
\eq{
 &\bar{BR}_t^{m,j}(\underline{\pi_t},\gamma_t^{l},\gamma_t^{m,-j}) :=\big\{ \tgamma_t^{m,j}: \forall x_t^{m,j}\in \cX^{m,j}, \tgamma_t^{m,j}(\cdot|x_t^{m,j})\nn\\
 &\in  \arg\max_{\gamma^{m,j}_t(\cdot|x_t^{m,j})}\E^{\gamma^{m,j}_t(\cdot|x_t^{m,j}) {\gamma}^{l}_t,\gamma_t^{m,-j},\,\underline{\pi_t}} \nn\\
 & 
\big\{ R_t^{m,j}( X^{l,m}_t,A^{l,m}_t) +\delta V_{t+1}^{m,j}(\underline{F}(\underline{\pi_t},\gamma_t^{l},\tgamma_t^{m,j},\gamma_t^{m,-j},A^{l,m}_t), X^{m,j}_{t+1}) \big\lvert \underline{\pi_t}, x_t^{m,j} \big\}  \big\}, \label{eq:m_P_mj2}
}

and 
\eq{
\pi_t(\cdot)=P^{\tsigma_{1:t-1}^l,\hat{\sigma}_{1:t-1}^m}(\cdot|a_{1:t-1}). \label{eq:pi_def}
}

\end{lemma}
\begin{proof}
\seq{
\eq{
& \E^{\gamma_t^{l},\hat{\gamma}_t^{m},\,\pi_t} \big\{ R_t^{m,j}(X^{l,m}_t,A^{l,m}_t) + V_{t+1}^{m,j}(\underline{F}(\underline{\pi_t},\gamma_t^{l},\hat{\gamma}_t^m,A^{l,m}_t), X^{m,j}_{t+1}) \big\lvert \underline{\pi_t}, x_t^{m,j} \big\}  \big\}
}
\eq{
&=\max_{\gamma^{m,j}_t(\cdot|x_t^{m,j})} \E^{\gamma_t^{l},\gamma_t^{m,j}(\cdot|x_t^{m,j}),\hat{\gamma}_t^{m,-j} \,\pi_t} \big\{ R_t^{m,j}(X^{l,m}_t,A^{l,m}_t) +\nn\\
&V_{t+1}^{m,j}(\underline{F}(\underline{\pi_t},\gamma_t^{l},\hat{\gamma}_t^m,A^{l,m}_t), X^{m,j}_{t+1}) \big\lvert \underline{\pi_t}, x_t^{m,j} \big\}  \big\}\nn\\
%
%
 %
&= \max_{\gamma^{m,j}_t(\cdot|x_t^{m,j})} \E^{ \gamma_t^{l},\gamma_t^{m,j}(\cdot|x_t^{m,j})\hat{\gamma}_t^{m,-j},\pi_t} \big\{ R_t^{m,j}(X^{l,m}_t,A^{l,m}_t)+ \E^{\tsigma_{t+1:T}^{l,m},\underline{F}(\underline{\pi}_t, \gamma_t^{l},\hat{\gamma}_t^m, A^{l,m}_t)}\nn\\
& \big\{\sum_{n=t+1}^T \delta^{n-t}R_n^{m,j}(X^{l,m}_n,A^{l,m}_n) |(\underline{\pi_t}, \gamma_t^{l},\hat{\gamma}_t^m),a^{l,m}_{1:t-1},A^{l,m}_t,x_{1:t}^{m,j},X_{t+1}^{m,j}\big\}\big\lvert \pi_t,x_t^{m,j}\big\}\label{b_eq:P1}\\
&= \max_{\gamma^{m,j}_t(\cdot|x_t^{m,j})} \E^{ \gamma_t^{l},\gamma_t^{m,j}(\cdot|x_t^{m,j})\hat{\gamma}_t^{m,-j}) \,\pi_t} \big\{ R_t^{m,j}(X^{l,m}_t,A^{l,m}_t)+\nn\\
&\max_{\sigma_{t+1:T}^{m,j}}\E^{\tsigma_{t+1:T}^{l},\sigma_{t+1:T}^{m,j},\tsigma_{t+1:T}^{m,-j},\underline{F}(\underline{\pi}_t, \gamma_t^{l},\hat{\gamma}_t^m, A^{l,m}_t)}\nn\\
& \big\{\sum_{n=t+1}^T \delta^{n-t}R_n^{m,j}(X^{l,m}_n,A^{l,m}_n) |(\underline{\pi_t}, \gamma_t^{l},\hat{\gamma}_t^m),a^{l,m}_{1:t-1},A^{l,m}_t,x_{1:t}^{m,j},X_{t+1}^{m,j}\big\}\big\lvert \pi_t,x_t^{m,j}\big\}\label{b_eq:P2}\\
%
%
%
&=\max_{\gamma^{m,j}_t(\cdot|x_t^{m,j})} \E^{\gamma_t^{l},\gamma_t^{m,j}(\cdot|x_t^{m,j})\hat{\gamma}_t^{m,-j},\,\pi_t} \big\{ R_t^{m,j}(X^{l,m}_t,A^{l,m}_t) +\nn\\ &\max_{\sigma_{t+1:T}^{m,j}}\E^{\sigma_t^{l},\hat{\sigma}_t^m,\tsigma_{t+1:T}^{l}\sigma_{t+1:T}^{m,j},\sigma_{t+1:T}^{m,-j},\underline{\pi}_{t}}  \nn\\
&
\big\{ \sum_{n=t+1}^T R_n^{m,j}(X^{l,m}_n,A^{l,m}_n) \big\lvert (\underline{\pi_t}, \gamma_t^{l},\hat{\gamma}_t^m),a^{l,m}_{1:t-1},A^{l,m}_t,x_{1:t}^{m,j},X_{t+1}^{m,j}\big\} \big\vert \pi_t,a^{l,m}_{1:t-1}, x_{1:t}^{m,j} \big\}\label{b_eq:P3} \\
&=\max_{\sigma_t^{m,j}} \E^{\sigma_t^{l},\sigma_t^{m,j}\hat{\sigma}_t^{m,-j},\underline{\pi}_{t}} \big\{ R_t^{m,j}(X^{l,m}_t,A^{l,m}_t) +\nn\\
&\max_{\sigma_{t+1:T}^{m,j}}\E^{\sigma_t^{l},\hat{\sigma}_t^m,\tsigma_{t+1:T}^{l}\sigma_{t+1:T}^{m,j},\sigma_{t+1:T}^{m,-j},\underline{\pi}_{t}} \nn\\
&\big\{ \sum_{n=t+1}^T R_n^{m,j}(X^{l,m}_n,A^{l,m}_n) \big\lvert F(\underline{\pi_t}, \gamma_t^{l},\hat{\gamma}_t^m),a^{l,m}_{1:t-1},A^{l,m}_t,x_{1:t}^{m,j},X_{t+1}^{m,j}\big\} \big\vert\pi_t, a^{l,m}_{1:t-1}, x_{1:t}^{m,j} \big\}\label{b_eq:P4} \\
&=\max_{\sigma_{t:T}^{m,j}} \E^{\sigma_t^{l},\sigma_t^{m,j}\hat{\sigma}_t^{m,-j},\tsigma_{t+1:T}^{l}\sigma_{t+1:T}^{m,j},\tsigma_{t+1:T}^{m,-j},\underline{\pi}_{t}}\big\{ \sum_{n=t}^T \delta^{t-n} R_n^{m,j}(X^{l,m}_n,A^{l,m}_n) \big\vert \pi_t, a^{l,m}_{1:t-1}, x_{1:t}^{m,j} \big\}\label{b_eq:P5}
}
}
where \eqref{b_eq:P1} follows from Lemma~\ref{b_lemma:1} in Appendix~\ref{b_app:lemmas}, \eqref{b_eq:P3} follows from Lemma~\ref{b_lemma:3}, \eqref{b_eq:P4} follows from the fact that $\hat{\sigma}^{m}_t$ are of type $m$ and definition of $\hat{\gamma}_t^m$.
  This proves the theorem.

\end{proof}

\section{Lemmas for converse}
\label{app:Proof_Exist}

\begin{proof}
We prove this by contradiction. This implies there exists $\underline{\pi_t}$ such that either (a)~\eqref{eq:P1b} doesn't have a solution, or (b)~\eqref{eq:P2} doesn't have a solution. 

\bit{
\item[(a)]If~\eqref{eq:P1b} doesn't have a solution (concerning follower):
Suppose for any equilibrium generating function $\theta$ that generates $(\tsigma^{l},\tsigma^{m})$ through forward recursion, there exists $t\in\cT,  a_{1:t-1}^{l,m}$ such that for $\underline{\pi_t}(\cdot)= P^{\tsigma^{l,m,}}(\cdot|a_{1:t-1}^{l,m})$, \eqref{eq:P1b} is not satisfied for $\theta$
i.e. for $ \tgamma^{l}_t = \theta^{l}[\underline{\pi_t}] = \tsigma_t^{l}(\cdot|\underline{\pi_t},\cdot), \tgamma^{m}_t = \theta^{m}[\underline{\pi_t}] = \tsigma_t^{m}(\cdot|\underline{\pi_t},\cdot)$, $\exists j, x_t^{m,j}$ such that
\seq{
\eq{
\tgamma_t^{m,j}(\cdot|x_t^{m,j})\notin  \arg\max_{\gamma^{m,j}_t(\cdot|x_t^{m,j})}\E^{ {\tgamma}^{l}_t,\gamma^{m,j}_t(\cdot|x_t^{m,j})\tgamma_t^{m,-j},\,\underline{\pi_t}}\nn\\
 &\hspace{-7cm}  
\big\{ R_t^{m,j}( X_t^{l,m},A^{l,m}_t) +\delta V_{t+1}^{m,j}(\underline{F}(\underline{\pi_t},\tgamma_t^{l,m},A_t^{l,m}), X^{m,j}_{t+1}) \big\lvert \underline{\pi_t}, x_t^{m,j} \big\}  
  }
  Let $t$ be the first instance in the backward recursion when this happens. This implies $\exists\ \hat{\gamma}_t^{m,j}$ such that
  \eq{
  &\E^{{\tgamma}^{l}_t,\hat{\gamma}_t^{m,j}(\cdot|x_t^{m,j}),\tgamma_t^{m,-j},\,\underline{\pi_t}} 
\big\{ R_t^{m,j}( X^{l,m}_t,A_t^{l,m}) +\nn\\
&\delta V_{t+1}^{m,j}(\underline{F}(\underline{\pi_t},\tgamma_t^{l,m},A^{l,m}_t),X^{m,j}_{t+1}) \big\lvert \underline{\pi_t}, x_t^{m,j} \big\}  
  \nn\\
  &>  \E^{ {\tgamma}^{l}_t,\tgamma^{m,j}_t(\cdot|x_t^{m,j}), \tgamma^{m,-j}_t,\,\underline{\pi_t}} 
\big\{ R_t^{m,j}( X^{l,m}_t,A^{l,m}_t) +\nn\\
&\delta V_{t+1}^{m,j}(\underline{F}(\underline{\pi_t},\tgamma_t^{l,m},A^{l,m}_t),X^{m,j}_{t+1}) \big\lvert \underline{\pi_t}, x_t^{m,j} \big\}   \label{b_eq:E1}
  }
  This implies for $\hat{\sigma}^{m,j}_t(\cdot|a_{1:t-1}^{l,m},x_{1:t-1}^{m,j},\cdot) = \hat{\gamma}_t^{m,j}$,
  \eq{
  &\E^{\tsigma_{t:T}^{l,m},\underline{\pi_t}} \big\{ \sum_{n=t}^T R_n^{m,j}(X^{l,m}_n,A^{l,m}_n) \big\lvert \underline{\pi_t}, a^{l,m}_{1:t-1}, x_{1:t}^{m,j} \big\}
  \nn\\
  &= \E^{ \tsigma_t^{l,m}, \underline{\pi_t}} \big\{ R_n^{m,j}(X^{l,m}_t,A^{l,m}_t) + \E^{ \tsigma_{t:T}^{l,m},\underline{\pi_t}}\nn\\
  &\big\{ \sum_{n=t+1}^T R_n^{m,j}(X^{l,m}_n,A^{l,m}_n) \big\lvert \underline{\pi_t},a^{l,m}_{1:t-1},A^{l,m}_t, x_{1:t}^{m,j},X_{t+1}^{m,j} \big\}  \big\vert \underline{\pi_t},  a_{1:t-1}^{l,m}, x_{1:t}^{m,j} \big\}
\\
  &= \E^{ \tsigma_t^{l,m}, \,\underline{\pi_t}} \big\{ R_n^{m,j}(X^{l,m}_t,A^{l,m}_t) + \E^{ \tsigma_{t+1:T}^{l,m},\underline{F}(\underline{\pi_t},\tgamma_t^{l,m},A^{l,m}_t)}\nn\\
  &\big\{ \sum_{n=t+1}^T R_n^{m,j}(X^{l,m}_n,A^{l,m}_n) \big\lvert a_{1:t-1}^{l,m},A^{l,m}_t, x_{1:t}^{m,j},X_{t+1}^{m,j} \big\}  \big\vert \underline{\pi_t},  a_{1:t-1}^{l,m},x_{1:t}^{m,j} \big\} \label{b_eq:E2}
}
\eq{
  &=\E^{\tilde{\gamma}^{l}_t,\tgamma^{m,j}_t(\cdot|x_t^{m,j}),\tgamma^{m,-j}_t   \, \underline{\pi_t}} \big\{ R_n^{m,j}(X^{l,m}_t,A^{l,m}_t) +\nn\\
  &V_{t+1}^{m,j} (\underline{F}(\underline{\pi_t}, \tilde{\gamma}^{l,m}_t, A^{l,m}_t), X_{t+1}^{m,j}) \big\lvert\underline{\pi_t},x_t^{m,j} \big\} \label{b_eq:E3}
  \\
  &< \E^{\tsigma_t^l,\hat{\gamma}^{m,j}_t(\cdot|x_t^{m,j}),{\tsigma}_t^{m,-j}, \, \underline{\pi_t}} \big\{ R_n^{m,j}(X^{l,m}_t,A^{l,m}_t) + \nn\\
  &V_{t+1}^{m,j} (\underline{F}(\underline{\pi_t},\tilde{\gamma}^{l,m}_t, A^{l,m}_t),X_{t+1}^{m,j}) \big\lvert \underline{\pi_t}, x_t^{m,j} \big\}\label{b_eq:E4}
  \\
  &= \E^{ \tsigma_t^{l},\hat{\sigma}_t^{m,j},{\tsigma}_t^{m,-j},\underline{\pi_t}} \big\{ R_n^{m,j}(X^{l,m}_t,A^{l,m}_t) +\E^{\tsigma_{t+1:T}^{l,m} \underline{F}(\underline{\pi_t},\tgamma_t^{l,m},A^{l,m}_t)}\nn\\
  &\big\{ \sum_{n=t+1}^T R_n^{m,j}(X^{l,m}_n,A^{l,m}_n) \big\lvert  a^{l,m}_{1:t-1},A^{l,m}_t, x_{1:t}^{m,j},X_{t+1}^{m,j}\big\} \big\vert \underline{\pi_t},a_{1:t-1}^{l,m}, x_{1:t}^{m,j} \big\}\label{b_eq:E5}
  \\
  &=\E^{ \tsigma_{t:T}^{l},\hat{\sigma}_t^{m,j},\tsigma_t^{m,-j},\underline{\pi_t}} \big\{ \sum_{n=t}^T R_n^{m,j}(X^{l,m}_n,A^{l,m}_n) \big\lvert \underline{\pi_t}, a_{1:t-1}^{l,m}, x_{1:t}^{m,j} \big\},\label{b_eq:E6}
  }
  where \eqref{b_eq:E2} follows from Lemma~\ref{b_lemma:3}, \eqref{b_eq:E3} follows from the definitions of $\tgamma_t^{m,j}$ and $\underline{\pi_t}$ and Lemma~\ref{b_lemma:1}, \eqref{b_eq:E4} follows from \eqref{b_eq:E1} and the definition of $\hat{\sigma}_t^{m,j}$, \eqref{b_eq:E5} follows from Lemma~\ref{b_lemma:2}, \eqref{b_eq:E6} follows from Lemma~\ref{b_lemma:3}. However, this leads to a contradiction since $(\tsigma^{l,m})$ is a SPE of the game.
}
\seq{
\item[(b)] ~\eqref{eq:P2}
 doesn't have a solution (concerning leader)

Suppose for any equilibrium generating function $\theta$ that generates $(\tsigma^{l,m})$ through forward recursion, there exists $t\in\cT,  a_{1:t-1}^{l,m}$ such that for $\underline{\pi_t}(\cdot)= P^{\tsigma^{l,m}}(\cdot|a_{1:t-1}^{l,m})$, \eqref{eq:P2} is not satisfied for $\theta$
i.e. for $ \tgamma^{l,m}_t = \theta^{l,m}[\underline{\pi_t}] = \tsigma_t^{l,m}(\cdot|\underline{\pi_t},\cdot)$, $\exists  x^{l}_t$ such that
\eq{
&\tgamma_t^{l} \notin \nn\\
&\arg\max_{\gamma_t^{l}(\cdot|x_t^l)}\E^{ {\gamma}^{l}_t,\bar{\gamma}_t^m,} \big\{ R_t^{l}(X^{l,m}_t,A^{l,m}_t)+\delta V_{t+1}^{l}(\underline{F}(\underline{\pi_t},{\tgamma}^{l}_t,\bar{\gamma}_t^m,A^{l,m}_t),X_{t+1}^{l})|\underline{\pi_t},x^{l}_t\big\},  \label{c_eq:P2}\\
&\text{where } \bar{\gamma}_t^m\in \bar{BR}_t^m(\underline{\pi_t}, \gamma_t^{l},\bar{\gamma}_t^m),\;\;\;
} 
  Let $t$ be the first instance in the backward recursion when this happens. This implies $\exists i\ \breve{\gamma}_t^{l}$ such that
  \eq{
  &\E^{ \breve{\gamma}^{l}_t,\hat{\gamma}_t^m,\,\underline{\pi}_t} \big\{ R_t^{l}(X^{l,m}_t,A^{l,m}_t)+\delta V_{t+1}^{l}(\underline{F}(\underline{\pi_t},{\tgamma}_t^{l},\hat{\gamma}_t^{m},A^{l,m}_t)|\underline{\pi_t},x^{l}_t\big\}  
  \nn\\
  &>  \E^{ {\tgamma}^{l,m}_t,\,\underline{\pi}_t} \big\{ R_t^{l}(X^{l,m}_t,A^{l,m}_t) +\delta V_{t+1}^{l}(\underline{F}(\underline{\pi_t},\tgamma^{l,m}_t,A^{l,m}_t),X_{t+1}^{l})|\underline{\pi_t},x^{l}_t\big\}   \label{c_eq:E1}
  }
  \eq{
 &\text{where }
 \hat{\gamma}_t^m\in \bar{BR}_t^m(\underline{\pi_t}, \breve{\gamma}^{l}_t,\hat{\gamma}_t^m)
 }
  This implies for $\breve{\sigma}^{l}_t(\cdot|a_{1:t-1}^{l,m},x_{1:t-1}^{l},\cdot) = \breve{\gamma}_t^{l}$,
  \eq{
  &\E^{ \tsigma_{t:T}^{l,m},\underline{\pi_t}} \big\{ \sum_{n=t}^T R_n^{l}(X^{l,m}_n,A^{l,m}_n) \big\lvert \underline{\pi_t}, a_{1:t-1}^{l,m}, x_{1:t}^{l} \big\}
  \nn\\
  &= \E^{ \tsigma_t^{l,m},\underline{\pi_t}} \big\{ R_t^{l}(X^{l,m}_t,A^{l,m}_t) + \nn\\
  &\E^{ \tsigma_{t:T}^{l,m},\underline{\pi_t}}\big\{ \sum_{n=t+1}^T R_n^{l}(X^{l,m}_n,A^{l,m}_n) \big\lvert \underline{\pi}_t,a^{l,m}_{1:t-1},A^{l,m}_t, x_{1:t}^{l},X_{t+1}^i \big\}  \big\vert \underline{\pi}_t,  a_{1:t-1}^{l,m}, x_{1:t}^{l} \big\} \label{a_eq:E2a}\\
   &= \E^{\tsigma_t^{l,m}, \,\underline{\pi_t}} \big\{ R_t^{l}(X^{l,m}_t,A^{l,m}_t) + \E^{ \tsigma_{t+1:T}^{l,m},\underline{F}(\underline{\pi_t},\tgamma_t^{l,m},A^{l,m}_t)}\nn\\
  &\big\{ \sum_{n=t+1}^T R_n^{l}(X^{l,m}_n,A^{l,m}_n) \big\lvert \underline{\pi_t},  , a^{l,m}_{1:t-1},A^{l,m}_t, x_{1:t}^{l},X_{t+1}^i \big\}  \big\vert \underline{\pi_t},  a_{1:t-1}^{l,m},x_{1:t}^{l} \big\} \label{c_eq:E2}
  \\
  &=\E^{ \tilde{\gamma}^{l,m}_t, \, \underline{\pi_t}} \big\{ R_t^{l}(X^{l,m}_t,A^{l,m}_t) + V_{t+1}^{l} (\underline{F}(\underline{\pi_t}, \tilde{\gamma}^{l,m}_t, A^{l,m}_t), X_{t+1}^{l}) \big\lvert\underline{\pi_t},x^{l}_t \big\} \label{c_eq:E3}
  \\
  &< \E^{ \breve{\gamma}^{l}_t,\hat{\gamma}^m_t, \underline{\pi_t}} \big\{ R_t^{l}(X^{l,m}_t,A^{l,m}_t) +V_{t+1}^{l} (\underline{F}(\underline{\pi_t}, {\tgamma}^{l}_t,\hat{\gamma}^m_t, A^{l,m}_t),X_{t+1}^{l}) \big\lvert \underline{\pi_t}, x^{l}_t \big\}\label{c_eq:E4}\\
  &= \E^{ \breve{\sigma}^{l}_t,\hat{\sigma}^m_t, \underline{\pi_t}} \big\{ R_t^{l}(X^{l,m}_t,A^{l,m}_t) +  \E^{ \tsigma_{t+1:T}^{l,m},F({\pi}_t, {\tgamma}^{l}_t,\hat{\gamma}^m_t, A^{l,m}_t)}\nn\\
  &\big\{ \sum_{n=t+1}^T R_n^{l}(X^{l,m}_n,A^{l,m}_n) \big\lvert a^{l,m}_{1:t-1},A^{l,m}_t, x_{1:t}^{l},X_{t+1}^{l}\big\} \big\vert \underline{\pi}_t,a_{1:t-1}^{l,m}, x_{1:t}^{l} \big\}\label{c_eq:E5}
}
\eq{
  &=\E^{\breve{\sigma}^{l}_t,\tsigma_{t+1:T}^{l},\hat{\sigma}^m_t,\tsigma_{t+1:T}^{m},\underline{\pi_t}} \big\{ \sum_{n=t}^T R_n^{l}(X^{l,m}_n,A^{l,m}_n) \big\lvert \underline{\pi_t}, a_{1:t-1}^{l,m}, x_{1:t}^{l} \big\},\label{c_eq:E6}
  }
  where \eqref{c_eq:E2} follows from~Lemma~\ref{l_lemma:3}, \eqref{c_eq:E3} follows Lemma~\ref{l_lemma:1}, \eqref{c_eq:E4} follows from \eqref{c_eq:E1}, \eqref{c_eq:E5} follows from Lemma~\ref{l_lemma:2}, \eqref{c_eq:E6} again follows from \\Lemma~\ref{l_lemma:3}. However, this leads to a contradiction since $(\tsigma^{l,m})$ is an SPE of the game.
  }

}
 
\end{proof}

\bibliographystyle{IEEEtran}

\end{document}